\documentclass[a4paper]{article}

\usepackage{amsfonts}
\usepackage{amssymb}
\newtheorem{st}      {Theorem}    [section]
\newtheorem{prop}[st]{Proposition}
\newtheorem{lem} [st]{Lemma}     
\newtheorem{cor} [st]{Corollary}

\newtheorem{exam}[st]{Example}   
\newtheorem{num}[st]{}

\newcommand{\ii}{\ensuremath{\mathcal{T}}}
\newcommand{\rr}{\ensuremath{\mathcal{R}}}
\newcommand{\thp}{\ensuremath{\theta}}

\newcommand{\f}{\ensuremath{\varphi}}
\newcommand{\rmap}{\longrightarrow}
\newcommand{\Boxe}{\raisebox{.8ex}{\framebox}}
\newcommand{\lmap}{\longleftarrow}

\newcommand{\U}{\ensuremath{\mathcal{U}}}
\newcommand{\OO}{\ensuremath{\mathcal{O}}}

\newcommand{\Ka}{\ensuremath{\mathcal{K}}}
\newcommand{\Ha}{\ensuremath{\mathcal{H}}}
\newcommand{\G}{\ensuremath{\mathcal{G}}}

\newcommand{\A}{\ensuremath{\mathcal{A}}}

\newcommand{\R}{\ensuremath{\mathcal{R}}}
\newcommand{\E}{\ensuremath{\mathcal{E}}}

\newcommand{\el}{\ensuremath{\mathcal{L}}}
\newcommand{\F}{\ensuremath{\mathcal{F}}}
\newcommand{\es}{\ensuremath{\mathcal{S}}}
\newcommand{\B}{\ensuremath{\mathcal{B}}}
\newcommand{\C}{\ensuremath{\mathcal{C}}}
\newcommand{\D}{\ensuremath{\mathcal{D}}}

\newcommand{\er}{\ensuremath{\mathcal{R}}}

\newcommand{\shG}{\ensuremath{\underline{Sh}(\G)}}
\newcommand{\abG}{\ensuremath{\underline{Ab}(\G)}}

\newcommand{\mod}{\ensuremath{\underline{Mod}}}
\newcommand{\shK}{\ensuremath{\underline{Sh}(\Ka)}}
\newcommand{\abK}{\ensuremath{\underline{Ab}(\Ka)}}

\newcommand{\cii}{\ensuremath{\C^{\infty}}}

\newcommand{\gc}{\ensuremath{\Gamma_c}}
\newcommand{\ps}{{\raise 1pt\hbox{\tiny (}}}

\newcommand{\pss}{{\raise 1pt\hbox{\tiny [}}}
\newcommand{\pdd}{{\raise 1pt\hbox{\tiny ]}}}
\newcommand{\pd}{{\raise 1pt\hbox{\tiny )}}}

\newcommand{\bs}{{\raise 1pt\hbox{\tiny [}}}
\newcommand{\bd}{{\raise 1pt\hbox{\tiny ]}}}
\newcommand{\nG}[1]{\ensuremath{\G^{\,\ps #1\,\pd}}}
\newcommand{\nK}[1]{\ensuremath{\Ka^{\,\ps #1\,\pd}}}

\def\cross{\mathinner{\mathrel{\raise0.8pt\hbox{$\scriptstyle>$}}
                 \joinrel\mathrel\triangleleft}}
\newcommand{\nH}[1]{\ensuremath{\Ha^{\,\ps #1\,\pd}}}
\newcommand{\nB}[1]{\ensuremath{B^{\ps #1\pd}}}

\newcommand{\Z}{\ensuremath{\mathcal{Z}}}

\def\cdot{{\raise 0.5pt\hbox{$\scriptscriptstyle\bullet$}}}
\def\compose{{\raise 1pt\hbox{$\scriptscriptstyle\circ$}}}
\def\dcross{{\raise 0.5pt\hbox{$\scriptscriptstyle\boxtime$}}}

\usepackage[all]{xy}
\usepackage{epsf} 
\usepackage{amsfonts}
\usepackage{amssymb}

\begin{document}

\title{A Homology Theory for \'Etale Groupoids\thanks{Research supported by NWO}}
\author { by Marius Crainic and Ieke Moerdijk}
\date{Preprint nr. 1605, May 1998, Utrecht University}
\pagestyle{myheadings}
\maketitle
\begin{abstract}
\'Etale groupoids arise naturally as models for leaf spaces of foliations, for orbifolds, and for orbit spaces of discrete group actions. In this paper we introduce a sheaf homology theory for \'etale groupoids. We prove its invariance under Morita equivalence, as well as Verdier duality between Haefliger cohomology and this homology. We also discuss the relation to the cyclic and Hochschild homologies of Connes' convolution algebra of the groupoid, and derive some spectral sequences which serve as a tool for the computation of these homologies.\\
\newline 
\hspace*{.2in}{\bf Keywords}: \'etale groupoids, homology, duality, spectral sequences, cyclic homology,   \linebreak foliations.
\end{abstract}

\hspace*{2.5in}{\bf Introduction}\\

\hspace*{.15in}In this paper we introduce a homology theory for \'etale groupoids. \'Etale groupoids serve as model for structures like leaf spaces of foliations, orbifolds, and orbit spaces of actions by discrete groups. In this sense, \'etale groupoids should be viewed as generalized spaces.\\
\hspace*{.3in} In the literature one finds, roughly speaking, two different approaches to the study of \'etale groupoids. One approach is based on the construction of the convolution algebras associated to an \'etale groupoid, in the spirit of Connes' non-commutative geometry (\cite{CoOp, Co3}), and involves the study of cyclic and Hochschild homology and cohomology of these algebras (\cite{BrNi, Co3}). The other approach uses methods of algebraic topology such as the construction of the classifying space of an \'etale groupoid and its (sheaf) cohomology groups (\cite{Bo, Ha2, Fourier}).\\
\hspace*{.3in} Our motivation in this paper is twofold. First, we want to give a more complete picture of the second approach, by constructing a suitable homology theory which complements the existing cohomology theory. Secondly, we use this homology theory as the main tool to relate the two approaches.\\
\hspace*{.3in} Let us be more explicit: In the second approach, one defines for any \'etale groupoid $\G$ natural cohomology groups with coefficients in an arbitrary $\G$-equivariant sheaf. These were introduced in a direct way by Haefliger (\cite{difcoh}). As explained in \cite{Fourier}, they can be viewed as a special instance of the Grothendieck theory of cohomology of sites (\cite{SGA}), and agree with the cohomology groups of the classifying space of $\G$ (\cite{emb}). Moreover, these cohomology groups are invariant under Morita equivalences of \'etale groupoids. (This invariance is of crucial importance, because the construction of the \'etale groupoid modelling the leaf space of a given foliation involves some choices which determine the groupoid only up to Morita equivalence.) We complete this picture by constructing a homology theory for \'etale groupoids, again invariant under Morita equivalence, which is dual (in the sense of Verdier duality) to the existing cohomology theory. Thus, one result of our work is the extension of ``the six operations of Grothendieck''(\cite{SGA}) from spaces to leaf spaces of foliations.\\\hspace*{.3in}   Our homology theory of the leaf space of a foliation reflects some geometric properties of the foliation. For example, by integration along the fibers (leaves) it is related to the leafwise cohomology theory studied by   Alvarez Lopez, Hector and others (see \cite{alv} and the references cited there). It also shows that the Ruelle-Sullivan current of a measured foliation (see \cite{CoOp}) lives in Hafliger's (closed) cohomology. The results in \cite{BrNi, Cra} (see also Proposition \ref{perrii}) imply that our homology is also the natural target for the (localized) Chern character. (We plan to describe some of these connections more explicitly in a future paper.)\\
\hspace*{.3in} The homology theory also plays a central role in explaining the relation between the sheaf theoretic and the convolution algebra approaches to \'etale groupoids, already referred to above. Indeed, the various cyclic homologies of \'etale groupoids can be shown to be isomorphic to the homology of certain associated \'etale groupoids; it extends the previous results of Burghelea, Connes, Feigin, Karoubi, Nistor, Tsygan. This connection explains several basic properties of the cyclic and periodic homology groups, and leads to explicit calculations (\cite{Cra}). The previous work on the Baum-Connes conjecture for discrete groups, or for proper actions of discrete groups on manifolds, suggest that this homology will play a role in the Baum-Connes conjecture for \'etale groupoids.\\
\hspace*{.3in}  From an algebraic point of view, our homology theory is an extension of the homology of groups, while from a topological point of view it extends compactly supported cohomology of spaces. In this context, we should emphasize that even in the simplest examples, the \'etale groupoids which model leaf spaces of foliations involve manifolds which are neither separated nor paracompact. Thus, an important technical ingredient of our work is a suitable extension of the notions related to compactly supported section of sheaves to non-separated (non-paracompact) manifolds. For example, as a special case of our results one obtains the Verdier (and Poincar\'e) duality for non-separated manifolds. Our notion of compactly supported sections is also used in the construction of the convolution algebra of a (non-separated) \'etale groupoid. We believe that this extension to non-separated spaces has a much wider use that the one in this paper, and we have tried to give an accessible presentation of it in the appendix. The results in the appendix also play a central role in the calculation concerning the cyclic homology of \'etale groupoids in \cite{Cra}, and make it possible to extend the results of \cite{BrNi} for separated groupoids to the non-separated case.\\
\hspace*{.3in} We conclude this introduction with a brief outline of the paper.\\
\hspace*{.3in} In the first section we review the basic definitions and examples related to \'etale groupoids, and in the second section we summarize the sheaf cohomology of \'etale groupoids. These two sections serve as background, and do not contain any new results. Readers familiar with this background should immediately go to section 3, and consult the earlier sections for notational conventions.\\
\hspace*{.3in} In section 3, we present the definition of our homology theory and mention some of its immediate properties.\\
\hspace*{.3in} In section 4, a covariant operation $\f_{\, !}$ for any map $\f$ between \'etale groupoids is introduced, which can intuitively be thought of as a kind of ``integration along the fiber'' at the level of derived categories. We then prove a Leray spectral sequence for this operation. This spectral sequence is extremely useful. For example, we will use it to prove the Morita invariance of homology. It also plays a crucial role in many calculations in \cite{Cra}.\\
\hspace*{.3in} In section 5, we prove that the operation $\el\f_{\, !}$ has a right adjoint $\f^{\, !}$ at the level of derived categories, thus establishing Verdier duality. The Poincar\'e duality between (Haefliger) cohomology and (our) homology of \'etale groupoids is an immediate consequence.\\
\hspace*{.3in} In section 6, we summarize the main aspects of the relation to cyclic homology. This section is based on \cite{Cra}, to which we refer the reader for detailed proofs and further calculations.\\
\hspace*{.3in} In an appendix, we show how to adapt the definition of the functor $\Gamma_{c}(X; \A)$ (assigning to a space $X$ and a sheaf $\A$ the group of compactly supported sections) in such a way that all the properties (as expressed in \cite{Borel}, say) can be proved without using Hausdorffness and paracompactness of the space $X$. This appendix can be read independently from the rest of the paper.

\tableofcontents

\section{\'Etale groupoids}

\hspace*{.3in} In this section we review the definition of topological groupoids, fix the notations, 
and mention some of the main examples. \\
\hspace*{.3in} Recall first that a {\it groupoid} $\G$ is a (small) category in which every arrow is invertible. We 
will write $\nG{0}$ and $\nG{1}$ for the set of objects and the set of arrows in $\G$, respectively, 
and denote the structure maps by:

\begin{equation}
\label{rel1}
 \xymatrix {
\nG{1} \times_{\nG{0}} \nG{1}  \ar[r]^-{m} & 
\nG{1} \ar[r]^-{i} &
\nG{1} \ar[r]_-{t}\ar@<1ex>[r]^-{s} &
\nG{0} \ar[r]^-{u} &
\nG{1} }\ , 
\end{equation}

Here $s$ and $t$ are the source an target, $m$ denotes composition ($m\ps g, h\pd = g \compose h$), 
$i$ is the inverse ($i\ps g \pd = g^{-1}$) and for any $x\in \nG{0}, u\ps x \pd = 1_{x}$ is the unit 
at $x$. We write $g: x\rmap y$ or $x \stackrel{g}{\rmap} y$ to indicate that $g \in \nG{1}$ is an 
arrow with $s\ps g \pd = x$ and $t\ps g \pd = y$.\par
A {\it topological groupoid} $G$ is similarly given by topological spaces $\nG{0}$ and $\nG{1}$ and by 
continuous structure maps as in (\ref{rel1}). For a smooth groupoid, $\nG{0}$ and $\nG{1}$ are smooth manifolds, 
and these structure maps are smooth; moreover, one requires $s$ and $t$ to be submersions, so that the 
fibered product $\nG{1} \times_{\nG{0}} \nG{1}$ in (\ref{rel1}) is also a manifold.

\begin{num}
\emph{
{\bf Definition.} A topological (smooth) groupoid $\G$ as above is called \'etale if the source map $s: \nG{1} \rmap \nG{0}$ is a local homeomorphism (local diffeomorphism). This implies that all other structure maps in (\ref{rel1}) are also local homeomorphisms (local diffeomorphisms).
}
\end{num}

\begin{num}
\label{germ}
\emph{
{\bf Germs.} Any arrow $g: x\rmap y$ in an \'etale groupoid induces a germ $\tilde{g}: (U, x)\tilde{\rmap} (V, y)$ from a neighborhood $U$ of $x$ in $\nG{0}$ to a neighborhood $V$ of $y$. Indeed, we can define $\tilde{g} = t \compose \sigma$, where $x\in U \subset \nG{0}$ is so small that $s: \nG{1}\rmap \nG{0}$ has a section $\sigma: U\rmap \nG{1}$ with $\sigma\ps x \pd = g$. If $U$ is so small that $t\, |_{\sigma(U)}$ is also a homeomorphism (resp. diffeomorphism), then $\tilde{g}: U\tilde{\rmap} V$ is also a homeomorphism (resp. diffeomorphism). We will also write $\tilde{g}$ for the germ at $x$ of this map $\tilde{g}: U\tilde{\rmap} V$. Note that $\tilde{1}_{\, x}$ is the identity germ, and that $\widetilde{\ps hg \pd} = \tilde{h} \tilde{g}$\, if $g: x\rmap y$ and $h: y\rmap z$.
}
\end{num}

\begin{num}
\label{ex1}
\emph{
{\bf Examples} of \'etale groupoids. (Note that in examples 3 and 4, the space $\nG{1}$ is in general not Hausdorff.)\\
 \hspace*{.3in} $1$. Any topological space (manifold) $X$ can be viewed as an \'etale groupoid $\underline{X}$, with identity arrows only ($\underline{X}^{\,\ps 0\,\pd} = X = \underline{X}^{\,\ps 1\,\pd}$, etc.). We will often simply denote this groupoid by $X$ again.\newline
\hspace*{.3in} $2$. If a (discrete) group $\Gamma$ acts from the right on a space $X$, one can form a groupoid $X\cross \Gamma$ with $(X\cross \Gamma)^{\,\ps 0\,\pd} = X$ and $(X\cross \Gamma)^{\,\ps 1\,\pd} = X \times \Gamma$, by taking as arrows $x\lmap y$ those  $\gamma \in \Gamma$ with $y = x\gamma $. This groupoid is called the translation groupoid of the action.\\
\hspace*{.3in} $3$. (\cite{Ha2, Bo}) The Haefliger groupoid $\Gamma^{q}$ has $\mathbb{R}^{q}$ for its space of objects. An arrow $x\rmap y$ in $\Gamma^{q}$ is a germ of a diffeomorphism $(\mathbb{R}^{q}, x)\rmap (\mathbb{R}^{q}, y)$. This groupoid and its classifying space $B\Gamma^{q}$ (cf. \ref{nerve} below) play a central role in foliation theory.\\
\hspace*{.3in} $4$. (see, for example, \cite{Wi, CoOp, Fourier}) For a foliation $(M, \F)$ of codimension $q$, its holonomy groupoid $Hol(M, \F)$ can be represented by an \'etale groupoid $Hol_{T}(M, \F)$, depending on the choice of a ``complete transversal'' $T$, i.e. a submanifold $T \subset M$ of dimension $q$ which is transversal to the leaves and which meets every leaf at least once. Two different such transversals $T$ and $T\, '$ give Morita equivalent (see \ref{morita} below) \'etale groupoids $Hol_{T}(M, \F)$ and $Hol_{T\, '}(M, \F)$.\\
\hspace*{.3in} $5$. Any orbifold gives rise to a smooth \'etale groupoid. These groupoids $\G$ coming from orbifolds have the special property that $\ps s, t \pd: \nG{1}\rmap \nG{0} \times \nG{0}$ is a proper map (see \cite{MP}). Groupoids with this property are called proper. For a proper groupoid, $\nG{1}$ is Hausdorff whenever $\nG{0}$ is.\\
\hspace*{.3in} $6$. Let $\G$ be an \'etale groupoid. A right \G-space is a space $X$ equipped with a map $p: X \rmap \nG{0}$ and an action $ X\times_{\nG{0}} \nG{1} \rmap X , \, \ps x, g\pd \mapsto x\, g$ satisfying the usual identities . If $X$ is a right $\G$-space, one can construct a groupoid $X \cross \G$, with $(X \cross \G)_{0}= X$ and $(X \cross \G)_{1}= X \times_{\nG{0}} \nG{1}$ :\linebreak an arrow $x \lmap y$ in $X \cross \G$ is an arrow $p\ps x\pd \stackrel{g}{\lmap} p\ps y\pd$ with $y= x\, g$. (A similar construction applies of course to left $\G$-spaces.)
}
\end{num}

\begin{num}
\emph{
{\bf Homomorphisms.} Let $\G$ and $\Ka$ be \'etale groupoids. A {\it homomorphism} $\f : \Ka \rmap \G$ is given by two continuous (or smooth) maps $\f_{0} : \nK{0} \rmap \nG{0}$ and $\f_{1} : \nK{1} \rmap \nG{1}$ which commute with all the structure maps in (\ref{rel1}) (i. e. : $\f_{0} s\ps g\pd= s\f_{1} \ps g\pd, \f_{1} \ps g \compose h\pd= \f_{1} \ps g\pd \compose \f_{1} \ps h\pd ,$ etc.)
}
\end{num}

\begin{num}
\label{morita}
\emph{
{\bf Morita equivalence.} A homomorphism $\f : \Ka \rmap \G$ is called a {\it Morita} (or weak, or essential) {\it equivalence} if:\\  
\hspace*{.3in} $1$. The map $s \pi_{2}: \nK{0} \times_{\nG{0}} \nG{1} \rmap \nG{0}$, defined on the space of pairs $\ps y, g\pd \in \nK{0} \times \nG{1}$ with $t\ps g\pd= \f \ps y\pd$, is an \'etale surjection; \\ 
\hspace*{.3in} $2$. The square:
} 
\end{num}
\[ \xymatrix {
\nK{1} \ar[r]^-{\f_{1}} \ar[d]_-{\ps s, t\pd} & \nG{1} \ar[d]^-{\ps s,t \pd}\\
\nK{0}\times \nK{0} \ar[r]^-{\f_{0} \times \f_{0}} & \nG{0}\times \nG{0} }\]
is a fibered product. \par
We often write $\f: \Ka \tilde{\rmap} \G$ to indicate that $\f$ is such a Morita equivalence. Two groupoids $\G$ and $\Ha$ are said to be {\it Morita equivalent} if there are Morita equivalences $\Ha \tilde{\lmap}   \Ka \tilde{\rmap} \G$. This is a transitive relation. One generally considers the category of \'etale groupoids obtained by formally inverting the Morita equivalences. In this category, an arrow $\Ha \rmap \G$ is represented by two homomorphisms, as in:
\[ \Ha \tilde{\lmap}   \Ka \rmap \G \]
(see \cite{Fourier} for more details).

\begin{num}
\label{bundles}
\emph{
{\bf Bundles.} (see, for example, \cite{Bo, Ha1, HiSk, Fourier, Mr, Mrr}) Let $B$ be a ``base'' space, and $\G$ an \'etale groupoid. A left $\G$-bundle over $B$ consists of a space $P$ , a map $\pi : P \rmap B$, and a left action of $\G$ on $P$ (see \ref{ex1}.6) which respects $\pi$ in the sense that $\pi \ps ge\pd= \pi\ps e\pd$. The action is called principal if the canonical map between fibered products:
\[ \nG{1}\times_{\nG{0}} P \rmap P\times_{B} P , \ \ \ps g, e\pd \mapsto \ps ge, e\pd \]
is a homeomorphism. \\
\hspace*{.3in} If $B= \nK{0}$ is the space of objects of another groupoid $\Ka$, the bundle $P$ is said to be $\Ka$-equivariant if $P$ is also equipped with a right $\Ka$-action, which commutes with the left action by $\G$: $\ps g e\pd  h= g \ps e h\pd$; in this case the maps $P \rmap \nK{0}$ and $P \rmap \nG{0}$ are denoted by $s_{P}$ (``source'') and $t_{P}$ (``target''), respectively. For instance, any homomorphism $\f :\Ka \rmap\G$ induces a $\Ka$-equivariant principal $\G$-bundle:
\[ P\ps\f \pd = \nK{0}\cross_{\nG{0}}\nG{1} \]
(the space considered also in \ref{morita}.1), with $s_{P\ps\f \pd}\ps y, g\pd = y, t_{P\ps\f \pd}\ps y, g\pd = s\ps g\pd$.
The isomorphism classes of $\Ka$-equivariant principal $\G$-bundles $P$ can be viewed as ``generalized'' or ``Hilsum-Skandalis'' morphisms:
\[ P: \Ka \rmap \G \]
\hspace*{.3in} The category so obtained is equivalent to the category obtained by inverting the Morita equivalences (see \ref{morita}). Thus, showing that a certain construction is invariant under Morita equivalence is the same as showing that it is functorial on generalized morphisms.
}
\end{num}

\begin{num}
\label{nerve}
\emph{
{\bf Nerve and classifying space.} For an \'etale groupoid $\G$, we write $\nG{n}$ for the space of composable strings of arrows in $\G$: 
\[ x_{0} \stackrel{g_{1}}{\longleftarrow} x_{1} \stackrel{g_{2}}{\longleftarrow} . \ . \ . \ \stackrel{g_{n}}{\longleftarrow} x_{n} \]
\hspace*{.3in} For $n= 0, 1,$ this agrees with the notation for the space of objects and arrows of $\G$, already introduced. The spaces $\nG{n}  \, (n \geq 0)$ together form a simplicial space:
\begin{equation}
\label{rel2}
 \xymatrix {
 \  .\ .\ .\ \ar@<2ex>[r]\ar@<1ex>[r]\ar[r]\ar@<-1ex>[r] & \, \nG{2} \ar@<1ex>[r]\ar[r]\ar@<-1ex>[r] & \, \nG{1} \ar[r]\ar@<-1ex>[r] & \, \nG{0} },
\end{equation}
with the face maps $d_{i}: \nG{n} \rmap \nG{n-1}$ defined in the usual way:
 \[ d_{i} \ps  \ g_{1},. . ., g_{n} \pd = \left \{ \begin{array}{ll} 
                                              \ps  \ g_{2}, . . . , g_{n} \pd & \mbox{if $i=0$} \\
                                              \ps  \ g_{1}, . . . , g_{i}g_{i+1}, . . . , g_{n} \pd & \mbox{if $1 \leq i \leq n-1$} \\
                                              \ps  \ g_{1}, ... , g_{n-1} \pd & \mbox{if $i=n$}
                                                        \end{array}
                                            \right. .\]
Its (thick \cite{Segal}) geometric realization is the classifying space of $\G$, denoted $B\G$. This space $B\G$ classifies homotopy classes of principal $\G-$bundles (\cite{BufLor, Fourier}). A Morita equivalence $\varphi : \Ha \tilde{\rmap} \G $ induces a weak homotopy equivalence $B\Ha \tilde{\rmap} B\G$.
}
\end{num}

\begin{num}
\label{assumptions}
\emph{
{\bf Overall assumptions.} It is important to observe that in many relevant examples, the space $\nG{1}$ of arrows of an \'etale groupoid $\G$ is not Hausdorff, (cf. 3, 4 in \ref{ex1}). However, for any space $X$ in this paper we do assume that $X$ has an open cover by subsets $U \subset X$ which are each paracompact, Hausdorff, locally compact, and of cohomological dimension bounded by a number $d$ (depending on $X$ but not on $U$). These assumptions hold for any (non-separated) manifold of dimension $d$, and in particular for each of the spaces $\nG{n}$ associated to a smooth \'etale groupoid.
}
\end{num}

\section{Sheaves and cohomology}

\hspace*{.3in} In this section we review the definition and main properties of the cohomology groups $H^{n}(\G;\A)$ of an \'etale groupoid $\G$ with coefficients in a $\G$-sheaf $\A$. These groups have been studied by Haefliger (\cite{difcoh},\cite{man}). They can also be viewed as cohomology groups of the topos of $\G$-sheaves (Grothendieck-Verdier) and were discussed from this point of view in \cite{Fourier}.

\begin{num}
\emph{
{\bf $\G$-sheaves.} Let $\G$ be an \'etale groupoid. A $\G$-sheaf is a sheaf $\es$ on the space $\nG{0}$, on which $\nG{1}$ acts continuously from the right. In other words, $\es$ is a right $\G$-space (\ref{ex1}.6) for which the map $\es \rmap \G$ is \'etale (a local homeomorphism). A morphism of $\G$-sheaves  $\es \rmap \es^{'}$ is a morphism of sheaves which commutes with the action. We will write $Sh(\G)$ for the category of all $\G$-sheaves of sets, and $\abG$ for the category of abelian $\G$-sheaves. These categories have convenient exactness properties: it is well known that $\shG$ is a topos, and (hence) that $\abG$ is an abelian category with enough injectives. If $R$ is a ring, we write $\mod_{R}(\G)$ for the category of $\G$-sheaves of $\G$-modules. Thus $\abG= \mod_{\mathbb{Z}}(\G)$. Later, we will mostly work with the category $\mod_{\mathbb{R}}(\G)$ of $\G$-sheaves of real vector spaces.
}
\end{num}

\begin{num}
\label{ex2}
\emph{
{\bf Examples.}\\
\hspace*{.3in} $1.$ For any set or abelian group $A$ the corresponding constant sheaf on $\nG{0}$ can be equipped with the trivial $\G$-action. We will refer to $\G$-sheaves of this form as constant $\G$-sheaves; they are simply denoted by $A$ again.\\
\hspace*{.3in} $2.$ The sheaf $\A= \C_{\nG{0}}$ of germs of continuous real-valued functions on $\nG{0}$ has the natural structure of a $\G$-sheaf: if $g: x\rmap y$ in $\nG{1}$ and $\alpha \in \A_{y}$ is a germ at $g$, then $\alpha\cdot g$ is defined as the composition $\alpha \compose \tilde{g}$ (cf. \ref{germ}). Similarly if $\G$ is a smooth \'etale groupoid, the sheaf $\Omega^{n}_{\nG{0}}$ of differential n-forms on $\nG{0}$ has a structure of a $\G$-sheaf ($n \geq 0$).\\
\hspace*{.3in} $3.$ Let $E$ be a sheaf on $\nG{0}$ (no action). To $E$ we can associate a $\G$-sheaf $E[\G]= E \times_{\nG{0}} \nG{1}= \{ \ps e, g\pd : g: x \rmap y, e \in E_{y}\}$. The sheaf projection is the map $E[\G] \rmap \nG{0}$ given by $\ps e, g\pd \mapsto s\ps g\pd$, while the $\G$-action is given by composition, $\ps e, g\pd \cdot h= \ps e, g\compose h\pd$. Sheaves (isomorphic to ones) of this form are said to be free $\G$-sheaves. The freeness is expressed by the adjunction property:
\[ Hom_{\G}(E[\G], \es) = Hom_{\nG{0}}(E, \es) \]
for any $\G$-sheaf $\es$.\\
\hspace*{.3in} $4.$ Each of the spaces $\nG{n}$ in the nerve of $\G$ (cf. \ref{nerve}) has the structure of a $\G$-sheaf, with sheaf projection:
\[ \varepsilon_{n}: \nG{n} \rmap \nG{0}, \ \ (x_{0} \stackrel{g_{1}}{\longleftarrow} x_{1} \stackrel{g_{2}}{\longleftarrow} . \ . \ . \ \stackrel{g_{n}}{\longleftarrow} x_{n}) \mapsto x_{n} \ ,\]
and the $\G$-action given by composition, $\ps g_{1}, . . . , g_{n}\pd \cdot h= \ps g_{1}, . . . , g_{n}h\pd$. This $\G$-sheaf is denoted $F_{n-1}(\G)$. For $n \geq 1$ these sheaves are free, in fact $\G_{n+1}= \nG{n}[\G]$. The system of $\G$-sheaves: 
\begin{equation}
\label{rel2.1}
 \xymatrix {
 \  .\ .\ .\ \ar@<2ex>[r]\ar@<1ex>[r]\ar[r]\ar@<-1ex>[r] & \, F_{2}(\G) \ar@<1ex>[r]\ar[r]\ar@<-1ex>[r] & \,  F_{1}(\G) \ar[r]\ar@<-1ex>[r] & \,  F_{0}(\G)} ,
\end{equation}
has the structure of a simplicial $\G$-sheaf, whose stalk at $x \in \nG{0}$ is the nerve of the comma category $x/\G$. This stalk is a contractible simplicial set.\\
\hspace*{.3in} $5.$ For any $\G$-sheaf of sets $\es$, one can form the free abelian $\G$-sheaf $\mathbb{Z}\es$ ; the stalk of $\mathbb{Z}\es$ at $x \in \nG{0}$ is the free abelian group on the stalk $\es_{x}$. In particular, from (\ref{rel2.1}) we obtain a resolution: 
\begin{equation}
\label{rel2.free}
 .\ .\ .\ \stackrel{\delta}{\rmap} \mathbb{Z}F_{1}(\G) \stackrel{\delta}{\rmap} \mathbb{Z}F_{0}(\G) \rmap \mathbb{Z} \rmap 0 
\end{equation}
of the constant $\G$-sheaf $\mathbb{Z}$, where $\delta$ is defined by the alternating sums of the face maps in (\ref{rel2.1}).\\
\hspace*{.3in} $6.$ If $\nG{0}$ is a topological manifold of dimension $d$, recall that its orientation sheaf $or$ is given by $or(U)=H_{c}^{d}(U;\mathbb{R})^{\vee}$, (see e.g. \cite{BoTu, Iv}, and the Appendix for compactly supported cohomology in the case where $\nG{0}$ is non-Hausdorff). It has a natural $\G$-action: for any arrow $g: x\rmap y$ in $\G$, let $U_{x}$ and $U_{y}$ be neighborhoods of $x$ and $y$, so small that $s:\nG{1}\rmap\nG{0}$ has a section $\sigma$ through $g$ with $t\compose\sigma: U_{x}\tilde{\rmap} U_{y}$. Then $t\compose\sigma$ induces a map $H_{c}^{d}(U_{x})\tilde{\rmap} H_{c}^{d}(U_{y})$, so also a map $H_{c}^{d}(U_{x})^{\vee}\tilde{\rmap} H_{c}^{d}(U_{y})^{\vee}$. Hence by taking germs, it gives an action $or_{y}\rmap or_{x}$.\\
\hspace*{.3in} Note that if $\nG{0}$ is oriented (i.e. as a sheaf on $\nG{0}$, $or$ is isomorphic to the constant sheaf $\mathbb{R}$), it is not necessarily constant as a $\G$-sheaf. When it is (i.e. when $\nG{0}$ is orientable and any arrow $g: x\rmap y$ gives an orientation-preserving germ $\tilde{g}$, cf. \ref{germ}) we say that $\G$ is {\it orientable}.
}
\end{num}

\begin{num}
\label{morph}
\emph{
{\bf Morphisms.} A morphism of \'etale groupoids $\f :\Ka \rmap \G$ induces an evident functor:
\[ \f^{*} : \shG \rmap \shK \]
by pullback (and similarly an exact functor $\f^{*} : \abG \rmap \abK$). This functor has a right adjoint:
\[ \f_{*}: \shG \rmap \shK\ .\]
For an $\Ka$-sheaf $\es$, the sheaf $\f_{*}(\es)$ on $\nG{0}$ is defined for any open set $U \subset \nG{0}$ by:
\[ \f_{*}(\es) (U)= Hom_{\Ka}( \f/U, \es) .\]
Here $\f/U= \{ \ps y, g\pd : y \in \nK{0}, g: \f\ps y\pd \rmap x, x\in U \}$, with $\Ka$-sheaf structure given by $( y, g) h = ( y, g\compose \f\ps h\pd )$. The $\G$-action on this sheaf $\f_{*}(\G)$ is defined as follows: for $\xi \in \f_{*}(\es)_{x}$ and $g: x\,' \rmap x$, let $U_{x}$ be a neighborhood of $x$ so that $\xi$ is represented by an element $\xi \in \f_{*}(\es)(U_{x})$, and let $U_{x^{'}}$ be so small that $s: \nG{1} \rmap \nG{0}$ has a section $\sigma : U_{x'} \rmap \nG{1}$ through $g$ with $t \compose \sigma(U_{x'}) \subset U_{x}$. Then define $\xi g\in \f_{*}(\es)_{x'}$ to be the element represented by the morphism:
\[ \tau : \f/U_{x'} \rmap \es \ , \ (y, f: \f\ps y\pd \rmap z) \mapsto \xi(\sigma\ps z\pd \compose f). \]
These adjoint functors $\f^{*}$ and $\f_{*}$ together constitute a topos morphism:
\[ \f : \shK \rmap \shG .\]
If $\f: \Ka \tilde{\rmap} \G$ is a Morita equivalence, then this morphism is an equivalence of categories $\shK = \shG$. In fact, topos morphisms $\shK \rmap \shG$ correspond exactly to generalized morphisms $\Ka \rmap \G$, or equivalently, to pairs of homomorphisms $\Ka \tilde{\lmap}   \Ha \rmap \G$ (cf. \ref{morita},\ref{bundles}).
}
\end{num}

\begin{num}
\emph{
{\bf Invariant sections.} Let $\es$ be a $\G$-sheaf. A section $\sigma: \nG{0} \rmap \es$ is called invariant if $\sigma\ps y\pd g= \sigma\ps x\pd$ for any arrow $g: x\rmap y$ in $\G$. We write:
\[ \Gamma_{inv}(\G, \es) \]
for the set of invariant sections; it is an abelian group if $\es$ is an abelian sheaf. (In fact $\Gamma_{inv}(\G, \es)= \f_{*}(\es)$ where $\f: \G \rmap 1$ is the morphism into the trivial groupoid).
}
\end{num}

\begin{num}
\label{cohom}
\emph{
{\bf Cohomology.} For an abelian $\G$-sheaf $\A$, the cohomology groups $H^{n}(\G; \A)$ are defined as the cohomology groups of the complex:
\[ \Gamma_{inv}(\G; \ii^{0}) \rmap \Gamma_{inv}(\G; \ii^{1}) \rmap  .\ .\ .\ \]
where $\A \rmap \ii^{0} \rmap \ii^{1} \rmap .\ .\ . $\ is any resolution of $\A$ by injective $\G$-sheaves. In other words:
\[ H^{n}(\G; -) = R^{n}\Gamma_{inv}(\G, -) \ .\]
(Thus, $H^{n}(\G; \A)$ is simply the cohomology of the topos $\shG$ with coefficients in $\A$.) It is obvious that a homomorphism $\f: \Ka \rmap \G$ induces homomorphisms in cohomology:
\[ \f^{*} : H^{n}(\G ; \A) \rmap H^{n}(\Ka ; \f^{*}\A) \ \ \ps n \geq 0\pd\  .\]
If $\f$ is a Morita equivalence, these are isomorphisms, since $\shG = \shK $.
}
\end{num}

\begin{num}
\emph{
{\bf Leray spectral sequence.} For any morphism $\f : \Ka \rmap \G$ and any $\Ka$-sheaf $\A$, there is a Leray spectral sequence
}
\end{num}
\[ E^{p, q}_{2} = H^{p}(\G ; R^{q}\f_{*} \A)  \Longrightarrow H^{p+q}(\Ka ; \A) \ .\]
(The $\G$-sheaf $R^{q}\f_{*} \A$ can be explicitly described as the sheaf associated to the presheaf $U \rmap H^{q}(\f/U; \A)$ where $\f/U$ is the groupoid associated to the (right) action of $\Ka$ on the space $(\f/U)_{0}$ used in \ref{morph} (cf. \ref{ex1}.6). See \cite{SGA}.

\begin{num}
\label{basic}
\emph{
{\bf Basic spectral sequence.} Let $\G$ be an \'etale groupoid, and let $\A$ be a $\G$-sheaf. By pull-back along $\varepsilon_{n}: \nG{n} \rmap \nG{0}$ (see \ref{ex2}.4), $\A$ induces a sheaf $\varepsilon_{n}^{*}(\A)$ on $\nG{n}$ which we often simply denote by $\A$ again. Consider for each $p$ and $q$ the sheaf cohomology $H^{q}(\nG{p}, \A)$ of the space $\nG{p}$. For a fixed $q$, these form a cosimplicial abelian group, and there is a basic spectral sequence:
\[ H^{p}H^{q}(\nG{\cdot}, \A) \Longrightarrow  H^{p+q}(\G; \A)  \ .\]
(It arises from the double complex $\Gamma(\nG{p}, \varepsilon_{p}^{*}\ii^{q})$ where $\A \rmap \ii^{\cdot}$ is an injective resolution.) \\
\hspace*{.3in} It follows that if $0\rmap \A\rmap \A^{0}\rmap \A^{1}\rmap .\ .\ .\ $ is any resolution by $\G$-sheaves $\A^{q}$  with the property that $\varepsilon_{p}^{*}(\A^{q})$ is an acyclic sheaf on $\nG{p}$, then $H^{*}(\G; \A)$ can be computed by the double complex 
\[ \Gamma ( \nG{p} ; \varepsilon_{p}^{*}( \A^{q} )) \ .\]
}
\end{num}

\begin{num}
\emph{
{\bf $\check{C}$ech spectral sequence.} An open set $U \subset \nG{0}$ is called saturated if for any arrow $g: x\rmap y$ in $\G$, one has $s\ps g\pd \in U$ iff $t\ps g\pd \in U$. For such a $U$ there is an evident ``full'' subgroupoid $\G|_{U} \subset \G$, with $U$ as space of objects. If $\U$ is an open cover of $\nG{0}$ by saturated opens, there is a spectral sequence: 
\[ \check{H}^{p}(\U; \Ha^{q}(\A)) \Longrightarrow H^{p+q}(\G ; \A) \]
where $\Ha^{q}(\A)$ is the presheaf $U \mapsto H^{q}(\G|_{U}; \A|_{U})$.
}
\end{num}

\begin{num}
\emph{
{\bf Hypercohomology.} For a cochain complex $\A^{\cdot}$ of abelian $\G$-sheaves the hypercohomology groups $\mathbb{H}^{\, \,n}(\G; \A^{\cdot})$ are defined in the usual way, as the cohomology groups of the double complex $\Gamma_{inv}(\G; \ii^{\cdot})$ where $\A^{\cdot} \rmap \ii^{\cdot}$ is a quasi-isomorphism into a cochain complex of injectives. (If $\A^{\cdot}$ is concentrated in degree $0$ one recovers the ordinary cohomology defined in \ref{cohom}). For each $q \in \mathbb{Z}$ denote by $\Ha^{q}(\A^{\cdot})$ the $q$-th cohomology $\G$-sheaf of $\A^{\cdot}$. If $\A^{\cdot}$ is bounded below, there is a spectral sequence for hypercohomology analogous to the one in \ref{basic}:
\[ H^{p}\mathbb{H}^{\, q}(\nG{\cdot} ; \A^{\cdot}) \Longrightarrow \mathbb{H}^{\, p+q}(\G ; \A^{\cdot}) \ .\]
}
\end{num}

\begin{num}
\label{ext}
\emph{
{\bf Ext functor.} Recall that for any $\G$-sheaf $\B$, the functor $Ext^{p}(\B; -)$ is defined as the $p$-th right derived functor of the functor $Hom_{\G}(\B; -)$. Thus $H^{p}(\G; \A)= Ext^{p}(\mathbb{Z}, \A)$. For later purposes we recall Yoneda's description of $Ext^{p}(\B, \A)$ as the group as equivalence classes of ``extensions'':
\[ 0 \lmap \B \lmap \E_{1} \lmap .\ .\ .\ \lmap \E_{n} \lmap \A \lmap 0 \]
(see e.g. \cite{MaL}). By composition of exact sequences, one defines a cap product:
\[ Ext^{q}(\C, \B) \otimes Ext^{p}(\B, \A) \stackrel{\compose}{\longrightarrow}  Ext^{p+q}(\C, \A) \ .\]
The same applies of course to the category $\mod_{\mathbb{R}}(\G)$ of $\G$-sheaves of real vector spaces. We use the notation $Ext_{\mathbb{R}}^{p}(\B, \A)$ here. Recall also that, over $\mathbb{R}$, the tensor product defines a functor $Ext_{\mathbb{R}}^{p}(\B, \A) \rmap Ext_{\mathbb{R}}^{p}(\C \otimes \B,\C \otimes \A)$. This gives an easy description of the cap product in cohomology:
\[ H^{q}(\G; \A) \otimes H^{p}(\G; \B) \stackrel{\cup}{\longrightarrow} H^{p+q}(\G; \A \otimes_{\mathbb{R}}\B ) \]
as:
\[ Ext^{q}(\mathbb{Z}, \A) \otimes Ext^{p}(\mathbb{Z}, \B) \rmap Ext^{q}(\mathbb{Z}, \A) \otimes Ext^{p}(\A, \A \otimes_{\mathbb{R}}\B) \stackrel{\compose}{\rmap} Ext^{p+q}(\mathbb{Z}, A \otimes_{\mathbb{R}}\B) \ .\] 
}
\end{num}

\begin{num}
\label{Inthom}
\emph{
{\bf Internal hom.} For two $\G$-sheaves $\A$ and $\B$ the sheaf $\underline{Hom}(\A, \B)$ on $\nG{0}$ carries a natural $\G$-action, hence gives a $\G$-sheaf $\underline{Hom}_{\, \G}(\A, \B)$ (or simply $\underline{Hom}(\A, \B)$ again). We recall that:
\[ \Gamma_{inv}(\G, \underline{Hom}_{\, \G}(\A, \B)) = Hom(\A, \B) \]
is the group of action preserving homomorphisms, i.e. morphisms in the category $\abG$. The derived functor of:
\[ \underline{Hom}(\A, -) : \abG \rmap \abG \]
will be denoted by $R^{p}\underline{Hom}(\A, -)$ or by $\underline{Ext^{p}}(\A, -)$. 
}
\end{num}

\section{Homology}

\hspace*{.3in} In this section we will introduce the homology groups $H_{n}(\G; \A)$ for any \'etale groupoid $\G$ and any $\G$-sheaf $\A$. Among the main properties to be proved will be the invariance of homology under Morita equivalence. \\
\hspace*{.3in} For any Hausdorff space $X$, the standard properties of the functor which assigns to a sheaf $\es$ its group of compactly supported sections $\Gamma_c(X; \es)$ are well known and can be found in any book on sheaf theory. In the appendix, we show how to extend this functor to the case where $X$ is not necessarily Hausdorff, while retaining all the standard properties. We emphasize that throughout this paper, $\Gamma_c$ will denote this extended functor.\\
\hspace*{.3in} Let us fix an \'etale groupoid $\G$. The spaces $\nG{0}$ and $\nG{1}$ (and hence the spaces $\nG{n}$ for $n \geq 0$) are assumed to satisfy the general conditions of \ref{assumptions}, but we will not assume that $\G$ is Hausdorff. We write $d=$ cdim$(\nG{0})$ for the cohomological dimension of $\nG{0}$. Thus, for any $n \geq 0$ and any Hausdorff open set $U \subset \nG{n}$, the (usual) cohomological dimension  of $U$ is at most $d$.

\begin{num}
\label{bar}
\emph{
{\bf Bar complex.} Let $\A$ be a $\G$-sheaf, and assume that $\A$ is c-soft as a sheaf on $\nG{0}$ (we will briefly say that $\A$ is a ``c-soft $\G$-sheaf''). For each $n \geq 0$, consider the sheaf $\A_{n}= \tau_{n}^{\, *}(\A)$ on $\nG{n}$ constructed by pull-back along $\tau_{n}: \nG{n} \rmap \nG{0}, \tau_{n}( x_{0} \lmap . . . \lmap x_{n}) = x_{0}$. It is again a c-soft sheaf because $\tau_{n}$ is \'etale. The groups $\Gamma_{c}(\nG{n}, \A_{n})$ of compactly supported sections, introduced in the Appendix, together form a simplicial abelian group:
\begin{equation}
\label{rel4.1}
 \xymatrix {
B_{\cdot}(\G; \A): \  .\ .\ .\ \ar@<2ex>[r]\ar@<1ex>[r]\ar[r]\ar@<-1ex>[r] &\gc (\nG{2};\A_{2}) \ar@<1ex>[r]\ar[r]\ar@<-1ex>[r] & \gc(\nG{1};\A_{1}) \ar[r]\ar@<-1ex>[r] & \gc (\nG{0};\A_{0}) } 
\end{equation}
with face maps:
\begin{equation}
\label{rel4.2}
 d_{i}: \gc (\nG{n};\A_{n})  \rmap \gc (\nG{n-1};\A_{n-1})  
\end{equation}
defined as follows. First, for the face map $d_{i}: \nG{n} \rmap \nG{n-1}$ (cf. \ref{nerve}) there is an evident map (isomorphism in fact)\ $\A_{n} \rmap d_{i}^{\, *}(\A_{n-1}) $, whose stalk at $\stackrel{\rightarrow}{g}=(x_{0} \stackrel{g_{1}}{\lmap} . . . \stackrel{g_{n}}{\lmap} x_{n})$ is the identity map for $i \neq 0$ and the action by $g_{0}:  (\A_{n})_{\stackrel{\rightarrow}{g}} =  \A_{x_{0}} \rmap \A_{x_{1}} = d_{0}^{\, *}(\A_{n-1})_{\stackrel{\rightarrow}{g}}$ if $i= 0$. The map $d_{i}$ in (\ref{rel4.2}) is now obtained from this by summation along the fibres (see \ref{lset}):
} 
\end{num}

\[ \xymatrix {
\gc (\nG{n}; \A_{n})   \ar[r]^-{d_{i}} \ar[d] &  \gc (\nG{n-1}; \A_{n-1}) \\
\gc (\nG{n}; d_{i}^{\, *}(\A_{n-1})) \ar[r]^-{\sim} & \gc (\nG{n-1}, (d_{i})_{\, !}d_{i}^{\, *}\A_{n-1}) \ar[u]_-{\gc(\nG{n-1}, \sum_{d_{i}})} }\]
The homology groups $H_{n}(\G; \A)$ are defined as the homology groups of the simplicial abelian groups (\ref{rel4.1}), or equivalently, as those of the associated chain complex given by the alternating sum $\delta= \sum\ps -1\pd^{i}d_{i}$.\\
Similarly, any bounded below chain complex $\es_{\cdot}$ of c-soft sheaves gives rise to a double complex:
\begin{equation}
\label{rel4.3}
B_{\cdot}(\G; \es_{\cdot})
\end{equation}
and we define $\mathbb{H}_{\, \, n}(\G; \es_{\cdot})$ to be the homology of the associated total complex.

\begin{num}
\label{lem1}
{\bf Lemma.} Any quasi-isomorphism $\es_{\cdot} \rmap \ii_{\cdot}$ between bounded below chain complexes of c-soft $\G$-sheaves induces an isomorphism
\[ \mathbb{H}_{\, n}(\G; \es_{\cdot}) \stackrel{\sim}{\rmap} \mathbb{H}_{\, n}(\G; \ii_{\cdot}) \ .\]
\end{num}

{\emph Proof:} The spectral sequence of the double complex (\ref{rel4.3}) takes the form $E_{p,q}^{2}= H_{p}H_{q}(\G_{\cdot}, \es_{\cdot}) \Longrightarrow H_{p+q}(\G; \es_{\cdot})$, where the $E_{p,q}^{1}$-term is the homology $\mathbb{H}_{\, p}(\nG{p}; \es^{\cdot})$ of the complex $\gc (\nG{p}; \es_{\cdot})$. The lemma thus follows from \ref{qi}.\ \ $\Boxe$

\begin{num}
\label{cres}
\emph{
{\bf c-soft resolutions.} Let $\A$ be an arbitrary $\G$-sheaf. There always exists a resolution:
\begin{equation}
\label{rel4.4}
 0 \rmap \A \rmap \es^{0} \rmap .\ .\ .\ \rmap \es^{d} \rmap 0 
\end{equation}
by c-soft $\G$-sheaves. For example, since the category of $\G$-sheaves has enough injectives, one can take any injective resolution $0 \rmap \A \rmap \ii^{\cdot}$ and take $\es^{\cdot}$ to be the truncation $\tau_{\leq d}(\ii^{\cdot})$ (softness of $\es^{d}$ then follows as in \cite{Borel}, p.55). Or, one can use for $\ii^{\cdot}$ the flabby Godement resolution of $\A$ on the space $\nG{0}$ with its natural $\G$-action, and truncate it. In  the case of a smooth \'etale groupoid and working over $\mathbb{R}$, one also has the standard resolution:
\[ 0 \rmap \A \rmap \A \otimes \Omega^{0} \rmap \A \otimes \Omega^{1} \rmap .\ .\ .\   \]
}
\end{num}
obtained from the $\G$-sheaves $\Omega^{\cdot}$ of differential forms on $\nG{0}$. (Note that the last two resolutions are functorial in $\A$.) \\
\hspace*{.3in} Any resolution (\ref{rel4.4}) maps into the truncated injective one. And, similarly, given two resolutions $0 \rmap \A \rmap \es^{\cdot}$ and $0 \rmap \A \rmap \ii^{\cdot}$, there is a resolution $\rr^{\cdot}$ (e.g. the truncated injective one) and a diagram:
\begin{equation}
\label{rel4.5}
\xymatrix {
\A \ar[d] \ar[r] & \ii^{\cdot} \ar@{.>}[d] \\ 
\es^{\cdot} \ar@{.>}[r] & \rr^{\cdot} }
\end{equation}
which commutes up to homotopy.

\begin{num}
\label{defhom}
\emph{
{\bf Definition of homology.} Let $\A$ be an arbitrary $\G$-sheaf, and let $0 \rmap \A \rmap \es^{0} \rmap \ .\ .\ .\ \es^{d} \rmap 0$ be a c-soft resolution. Then $\es^{-\cdot}$ is a bounded chain complex (non-zero in degrees between $-d$ and $0$) and we define $H_{n}(\G; \A)$ to be $\mathbb{H}_{\, n}(\G; \es^{-\cdot})$. By \ref{cres} ((\ref{rel4.5})) and lemma \ref{lem1}, this definition is independent of the choice of the resolution. Observe that
\[ H_{n}(\G; \A) = 0 \ \ \ \ {\rm for\  all}  \ \  n < -d \ .\]
These groups can be viewed also as compactly supported cohomology groups (see \ref{extreme}.3 and \ref{comcoh} below).
}
\end{num}

\begin{num}
\label{extreme}
\emph{
{\bf Extreme cases.} \\
\hspace*{.3in} $1.$ If $\nG{0}$ is a point, i.e. if $\G$ is a discrete group, then $H_{n}(\G; \A)$ is the usual group homology of $\G$.\\
\hspace*{.3in} $2.$ If $\G$ is a discrete groupoid, $\nG{\cdot}$ is a simplicial set, and $H_{n}(\G; \A)$ is the usual simplicial homology of $\nG{\cdot}$ with twisted coefficients.\\
\hspace*{.3in} $3.$ If $\G$ is a Hausdorff space $X$ (viewed as a ``trivial'' groupoid, \ref{ex1}.1) then $H_{n}(\G; \A) = H_{c}^{-n}(X; \A)$ is the usual cohomology with compact supports (although graded differently). So the spectral sequence occurring in the proof of lemma \ref{lem1} could be written as $H_{p}H_{c}^{-q}(\nG{\cdot}; \A)) \Longrightarrow H_{p+q}(\G; \A)$.
}
\end{num}

\begin{num}
\label{les}
\emph{
{\bf Long exact sequence.} Any short exact sequence: 
\[ 0 \rmap \A \rmap \B \rmap \C \rmap 0 \]
of $\G$-sheaves induces a long exact sequence in homology:
\[ .\ .\ .\ \rmap H_{n+1}(\G; \C) \rmap H_{n}(\G; \A) \rmap H_{n}(\G; \B) \rmap H_{n}(\G; \C) \rmap .\ .\ .\ \]
The proof is standard. (The truncated Godement resolutions give a short exact sequence of resolutions $0 \rmap \es^{\cdot}(\A) \rmap \es^{\cdot}(\B) \rmap \es^{\cdot}(\C) \rmap 0$.)
}
\end{num}

\begin{num}
\label{functor}
\emph{
{\bf Functoriality.} Compactly supported cohomology of spaces (\ref{extreme}.3) is covariant along local homeomorphisms and contravariant along proper maps. Analogous properties hold for homology of  \'etale groupoids. Consider a homomorphism $\f: \Ka\rmap \G$ between \'etale groupoids.\\
\hspace*{.3in} $1.$ Suppose that $\f$ is proper, in the sense that each $\f_{n}: \nK{n} \rmap \nG{n}$ is a proper map (cf. \ref{lspr}). Then for any $\G$-sheaf $\A$ one obtains homomorphisms:
\[ \gc(\nG{n}; \A_{n}) \rmap \gc (\nK{n}; \f^{*}(\A)_{n}) \]
by pullback, and hence a homomorphism:
\[ \f^{*}: H_{n}(\G; \A) \rmap H_{n}(\Ka; \f^{*}\A) \ .\]
In other words, homology is contravariant along proper maps.\\
\hspace*{.3in} $2.$ Suppose $\f$ is \'etale, in the sense that each $\f_{n}: \nK{n} \rmap \nG{n}$ is a local homeomorphism (it is not difficult to see that the assumption is only about $\f_0$). Let $\es$ be a c-soft $\G$-sheaf. For the sheaf $\es_{n}= \tau_{n}^{\, *}(\es)$ on $\nG{n}$  summation along the fibers defines a homomorphism: 
\[ (\f_{n})_{\, !}\tau_{n}^{\, *}(\f^{*}(\es)) = (\f_{n})_{\, !}\f_{n}^{*}(\tau_{n}^{\, *}\es) \rmap \tau_{n}^{\, *}(\es)\ ,\]
and hence a homomorphism :
\[ \gc (\nK{n}; \f^{*}(\es)_{n}) \rmap \gc (\nG{n}; \es_{n})\ .\]
These homomorphisms, for each $n \geq 0$, commute with the face operators (\ref{rel4.2}). Since the functor $\f^{*}$ is (always) exact and preserves c-softness (because $\f$ is \'etale), this gives for each $\G$-sheaf $\A$ a homomorphism:
\[ H_{n}(\Ka; \f^{*}\A) \rmap H_{n}(\G; \A) \ .\]
\hspace*{.3in} $3.$ Suppose that $\f$ is \'etale, and moreover suppose that for each $n$ the square:
}
\end{num}
\begin{equation}
\label{rel4.6}
\xymatrix {
 \nK{n} \ar[d]_-{\f_{n}} \ar[r]^-{\tau_{n}} & \nK{0} \ar[d]^{\f_{0}} \\
\nG{n} \ar[r]^-{\tau_{n}} & \nG{0}  } 
\end{equation}
is a pullback. (Morphisms of this kind are exactly the projections $X\cross\G \rmap \G$ associated to \'etale $\G$-spaces $X$ .) For such a $\f$, there is an exact functor:
\[ \f_{\, !}: \abK \rmap \abG \]
which preserves c-softness. (at the level of underlying sheaves, it is simply the functor $(\f_{0})_{\, !}: \underline{Ab}(\nK{0}) \rmap \underline{Ab}(\nG{0})$ of \ref{losr}). For any c-soft $\Ka$-sheaf $\B$, there is a natural isomorphism: 
\[ \gc (\Ka_{n}; \B_{n}) = \gc (\Ka_{n}; \tau_{n}^{\, *}\B) = \gc (\nG{n}; (\f_{n})_{\, !}\tau_{n}^{\, *}\B) = \gc (\nG{n}; \tau_{n}^{\, *}(\f_{0})_{\, !}\B) = \gc (\nG{n}; \f_{\, !}(\B_{n})) ,\]
for any $n \geq 0$. These yield an isomorphism
\[ H_{n}(\Ka; \B) = H_{n}(\G; \f_{\, !}\B), \]
for any $\Ka$-sheaf $\B$.\\
\hspace*{.3in} Note that even if $\f$ is not \'etale, a functor $\f_{\, !}$ can be defined in this way (but it is no longer exact). See also \ref{wequiv}.4.

\begin{prop} Let $\f, \psi: \Ka\rmap\G$ be two \'etale homomorphisms, $\A\in\abG$ and $\thp: \psi\rmap\f$ a continuous transformation of functors. Denote by $\thp^{*}: H_{*}(\Ka; \f^{*}\A)\rmap H_{*}(\Ka; \psi^{*}\A)$ the map induced by the sheaf map $\f^{*}\A \rmap \psi^{*}\A, a\mapsto a\thp$, and by $\f_{*}, \psi_{*}$ the maps induced by $\f, \psi$ in homology (cf. \ref{functor}.2). Then:
\[ \xymatrix {
 H_{*}(\Ka; \f^{*}\A) \ar[dd]^-{\thp^{*}} \ar[dr]^-{\f_{*}}\\
 & H_{*}(\G; \A) & , \ \ \f_{*}= \psi_{*}\thp^{*}\\
H_{*}(\Ka; \psi^{*}\A) \ar[ur]_-{\psi_{*}} } \]
Moreover, the construction of $\thp^{*}$ is functorial with respect to $\thp$.
\end{prop}

\emph{Proof:} We may assume that $\A$ is c-soft. Then a homotopy between the maps:
\[ B_{\cdot}(\Ka; \f^{*}\A) \rmap B_{\cdot}(\G; \A) \]
inducing $\f_{*}$ and $\psi_{*}\thp^{*}$ in homology is given by:
\[ H= \sum_{i= 0}^{n} \ps -1\pd^i H_{i}: B_{n}(\Ka; \f^{*}\A) \rmap B_{n+ 1}(\G; \A),\]
where the $H_i$'s are defined as follows. Consider:
\[ h_i: \nK{n} \rmap \nG{n+ 1}, \]
\[ h_i(k_1, .\ .\ .\ , k_n)= \left \{ \begin{array}{ll}
                (\thp\ps t\ps k_1\pd\pd, \psi\ps k_1\pd, .\ .\ .\ , \psi\ps k_n\pd) & \mbox{if $i=0$} \\
                (\f\ps k_1\pd, .\ .\ .\, \f\ps k_i\pd, \thp\ps s\ps k_i\pd\pd, \psi\ps k_{i+ 1}\pd, .\ .\ .\ , \psi\ps k_n\pd) & \mbox{if $1 \leq i \leq n$}
                                      \end{array}
                              \right. . \]
Using the obvious (identity) isomorphisms $h_{i}^{*}(\A_{n+ 1}) \cong (\f^{*}\A)_{n}$, and summation along the fiber of the (\'etale) $h_{i}$'s (see \ref{lset} in Appendix), we get the homomorphisms:
\[ H_i: B_{\cdot}(\Ka; \f^{*}\A) \rmap B_{\cdot + 1}(\G; \A) \ .\]
The naturality with respect to $\thp$ is obvious. \ \ \ $\Boxe$

\begin{num}
\label{hyperhom}
\emph{
{\bf Hyperhomology.} Consider any bounded below chain complex $\A_{\cdot}$ of $\G$-sheaves. Let $\A_{\cdot} \rmap \R_{\cdot}$ be a q.i. into a bounded below chain complex of c-soft $\G$-sheaves. (Such an $\R_{\cdot}$ can be constructed for example by considering a resolution $\A_{\cdot} \rmap \es^{0}_{\cdot} \rmap .\ .\ .\ \es^{d}_{\cdot} \rmap 0$ as in \ref{cres} and then taking the total complex of the double complex $\es_{q}^{-p} \ (p, q \in \mathbb{Z}, -d \leq p \leq 0)$. Define the hyperhomology $\mathbb{H}_{\, *}(\G; \A_{\cdot})$ to be the homology of the total complex associated to the double complex $B_{\cdot}(\G; \R_{\cdot})$. This definition of $\mathbb{H}_{\, *}(\G; \A_{\cdot})$ does not depend on the choice of the resolution $\R_{\cdot}$ (cf. lemma \ref{lem1}).
}
\end{num}

\begin{prop}
\label{prop1}
(Hyperhomology spectral sequence) Let $\A_{\cdot}$ be a bounded below chain complex of $\G$-sheaves as above, and consider for each $q\in \mathbb{Z}$ the homology $\G$-sheaf $\Ha_{q}(\A_{\cdot})$. There is a spectral sequence:
\[ H_{p}(\G; \Ha_{q}(\A_{\cdot})) \Longrightarrow \mathbb{H}_{\, p+q}(\G; \A_{\cdot})\ . \]
\end{prop}

{\emph Proof}: Consider the truncated Godement resolution $0 \rmap \A_{\cdot} \rmap \es^{0}_{\cdot} \rmap .\ .\ .\ \rmap \es^{d}_{\cdot} \rmap 0$. It has the property that for each $q$, it also yields c-soft resolutions of the cycles $\Z_{q}$, the boundaries $\B_{q}$ and the homology $\Ha_{q}(\A_{\cdot})$. Write $\C$ for the triple complex:
\[ \C_{p,q,r} = \gc (\nG{p}; \es^{-r}_{q})\ ,\]
and let $\D$ be the double complex:
\[ \D_{n,q} = \bigoplus_{p+q= n}\C_{p,q,r}\ .\]
The total complex of $\C$, and hence also that of $\D$ , compute $\mathbb{H}(\G; \A_{\cdot})$. Furthermore, by the property of the resolution just mentioned (and  the fact that $\gc(\nG{p}; -)$ preserves exact sequences of c-soft sheaves) we have for fixed $p$ and $r$ that:
\[ H_{q}(\C_{p,\cdot,r}) = \gc (\nG{p}; \Ha_{q}(\es^{-r}_{\cdot}))\ .\]
Hence, for a fixed $n$,
\[ H_{q}(\D_{n,\cdot}) = \bigoplus_{p+r=n} \gc(\nG{p}; \Ha_{q}(\es^{-r}_{\cdot}))\ .\]
But $\Ha^{q}(\A_{\cdot}) \rmap \Ha_{q}(\es^{0}_{\cdot}) \rmap \Ha_{q}(\es^{1}_{\cdot}) \rmap .\ .\ .\ $ is a resolution of $\Ha^{q}(\A_{\cdot})$, so for a fixed $q$ the double complex $\gc (\nG{\cdot}; \Ha_{q}(\es^{-\cdot}))$ computes $H_{*}(\G; \Ha_{q}(\A_{\cdot}))$. Thus:
\[ H_{n}H_{q}(\D_{\cdot,\cdot}) = H_{n}(\G; \Ha_{q}(\A_{\cdot}))\ ,\]
and the desired spectral sequence is simply the spectral sequence $H_{n}H_{q}(\D) \Longrightarrow H_{n+q}(Tot(\D))$ for the double complex $\D$.\ \ $\Boxe$

\begin{num}
\label{cap}
\emph{
{\bf Cap product.} For an \'etale groupoid $\G$, the $Ext$-groups (\ref{ext}) act on the homology by a cap product:
\begin{equation}
\label{rel4.7}
 H_{n}(\G; \B) \otimes Ext^{p}(\B, \A) \stackrel{\cap}{\rmap} H_{n-p}(\G; \A)\ .
\end{equation}
\hspace*{.3in} For example, for $p=1$ an element of $Ext^{1}(\B, \A)$ can be represented by an exact sequence $0 \lmap \B \lmap \E \lmap \A \lmap 0$, which yields a boundary map $H_{n}(\G; \B) \rmap H_{n-1}(\G; \A)$ for the long exact sequence of \ref{les}. For $p > 1$, the cap product can be constructed in the same way (by decomposing a longer extension $0 \lmap \B \lmap \E_{1} \lmap .\ .\ .\ \lmap \E_{n} \lmap \A \lmap 0$ into short exact sequences). \\
\hspace*{.3in} In particular, when working over $\mathbb{R}$, this yields a simple description of the cap product relating homology and cohomology of \'etale groupoids:
\[ H_{n}(\G; \B) \otimes H^{p}(\G; \A) \stackrel{\cap}{\rmap} H_{n-p}(\G; \B\otimes_{\mathbb{R}} \A)\ .\]
\hspace*{.3in} The cap product satisfies the usual ``projection formula'' for a morphism $\alpha: \C \rmap \A$. Explicitly, $\alpha$ induces $\alpha_{*}: H_{*}(\G; \C) \rmap H_{*}(\G; \A)$ and $\alpha_{*}: Ext^{p}(\B; \C) \rmap Ext^{p}(\B; \A)$, and we have for any $u \in H_{n}(\G; \B)$ and $\xi \in Ext^{p}(\B; \C)$ that:
\[ \alpha_{*}(u \cap\xi) = u \cap\alpha_{*}(\xi)\ .\]
(For $p= 1$ this is just the naturality of the exact sequence \ref{les}).
}
\end{num}

\begin{num}
\emph{
{\bf Remark.} The $d^{2}$ boundary of the hyperhomology spectral sequence \ref{prop1}:
\[ \alpha_{p,q}^{2}: H_{p}(\G; \Ha_{q}(\A_{\cdot})) \rmap H_{p-2}(\G; \Ha_{q+1}(\A_{\cdot})) \]
is given by the cap product with an element $u_{q}(\A_{\cdot}) \in Ext^{2}(\Ha_{q}(\A_{\cdot}), \Ha_{q+1}(\A_{\cdot})).$ Let $\pi_{q}: \Z_{q}(\A_{\cdot}) \rmap  \Ha_{q}(\A_{\cdot})$ be the quotient map from the sheaf of cycles $\Z_{q}(\A_{\cdot})$. Then the extension 
\[ 0 \lmap \Ha_{q}(\A_{\cdot}) \stackrel{\pi_{q}}{\lmap} \Z_{q}(\A_{\cdot}) \stackrel{d}{\lmap} \A_{q+1} \lmap \Z_{q+1}(\A_{\cdot}) \lmap 0 \]
defines an element $v \in Ext^{2}(\Ha_{q}(\A_{\cdot}), \Z_{q+1}(\A_{\cdot}))$, and $u_{q}(\A_{\cdot})$ is $(\pi_{q+1})_{*}(\A_{\cdot})$. This is immediate from the construction of the spectral sequence (proof of \ref{prop1}), and the general description of the boundaries of the spectral sequence induced by a double complex.
}
\end{num}

\begin{num}
\label{etalecate}
\emph{
{\bf Remark.} Recall that a topological category $\G$ is said to be \'etale if all its structure maps are local homeomorphisms. Thus, such a category is given by maps as in (\ref{rel1}), except for the absence of an inverse $i:\nG{1}\rmap\nG{1}$. The definitions and the results of this section hold equally well for the more general context of such \'etale categories, and for this reason we have tried to write the proofs in such a way that they apply verbatim to this general context. The same is true for the next section, provided one takes sufficient care to define Morita equivalence for categories in the appropriate way. (In this paper we will only use the homology for \'etale categories in Proposition \ref{relc3}.) }
\end{num}

\section{Leray spectral sequence; Morita invariance}

\hspace*{.3in} In this section we construct for each morphism $\f: \Ka \rmap \G$ between \'etale groupoids a functor $\f_{\, !}$ from c-soft $\Ka$-sheaves to c-soft $\G$-sheaves. We derive a Leray spectral sequence for this functor (\ref{spseq}), of which the invariance of homology under Morita equivalences will be an immediate consequence (\ref{morinv}).

\begin{num}
\label{comma}
\emph{
{\bf Comma groupoids of a homomorphism.} Let $\f: \Ka \rmap \G$ be a homomorphism of \'etale groupoids. For each point $x \in \nG{0}$ consider the ``comma groupoid'' $x/\f$, whose objects are the pairs $(y, g: x \rmap \f\ps y\pd)$ where $y\in \nK{0}$ and $g\in \nG{1}$. An arrow $k: \ps y, g\pd \rmap \ps y\, ', g\, '\pd$ in $x/\f$ is an arrow $k: y\rmap y\, '$ in $\Ka$ such that $\f\ps k\pd\compose g= g\, '$. When equipped with the obvious fibered product topology, $x/\f$ is again an \'etale groupoid. It should be viewed as the fiber of $\f$ above $x$; more exactly, there is a commutative diagram (see also \ref{fibprod}):
\[ \xymatrix {
 x/\f \ar[r]^-{\pi_{x}} \ar[d] & \Ka \ar[d]^-{\f} \\
 1 \ar[r] \ar[r]^-{x} & \G } \]
Note that an arrow $g: x \rmap x\, '$ in $\G$ induces a homomorphism:
\begin{equation}
\label{rel5.0}
 g^{*}: x\, '/\f \rmap x/\f 
\end{equation}
by composition. Thus the groupoids $x/\f$ together form a right $\G$-bundle of groupoids. (If $\f_{0}: \nK{0}\rmap\nG{0}$ is a local homeomorphism, then it is a $\G$-sheaf of groupoids.)\\
\hspace*{.3in} More generally, for any $A \subset \nG{0}$ the comma groupoid $A/\f$ is defined by:
\[ (A/\f)^{\ps i\pd} = \bigcup_{x\in A} (x/\f)^{\ps i\pd} \subset \nK{i} \times_{\nG{0}} \nG{1} \ \ , i \in \{ 0, 1\} \]
(with the induced topology). The nerve of $A/\f$ consists of the spaces:
\[ (A/\f)^{\ps n\pd} = \{ (y_{0} \stackrel{k_{1}}{\lmap} .\ .\ .\ \stackrel{k_{n}}{\lmap}y_{n}, \, \f\ps y_{n}\pd \ \stackrel{g}{\lmap} x): k_{i}\in \nK{1}, g\in \nG{1}, x\in A \} \ .\]
When $\f = id: \G \rmap \G$, these are simply denoted by $x/\G$, $A/\G$. Dually one defines the comma groupoids $\f/x, \f/A, \G/x, \G/A$ (consisting on arrows ``going into $x$'').
}
\end{num}

\begin{num}
\label{flosr}
\emph{
{\bf The functors $\f_{\, !},\ L_n\f_{\, !},\ \el\f_{\, !}$.} Let $\f: \Ka \rmap \G$ be as above, and let $\A$ be a  $\Ka$-sheaf. We define a simplicial $\G$-sheaf $B_{\cdot}(\f; \A )$ in analogy with the definition of the bar-complex \ref{bar}. On the spaces $\nK{n} \times_{\nG{0}} \nG{1}$ (which form the nerve of $\nG{0}/\f$, cf. \ref{comma}) of strings of the form:
\[ \f \ps y_{0}\pd \stackrel{\f \ps k_{1}\pd }{\longleftarrow} .\ .\ .\ \stackrel{\f \ps k_{n}\pd }{\longleftarrow} \f \ps y_{n}\pd \stackrel{g}{\lmap} x \]
we define the maps:
\[ \alpha_{n} : \nK{n} \times_{\nG{0}} \nG{1} \rmap \nK{0} , \ \ps k_{1},\ .\ .\ .\ k_{n}, g\pd \mapsto t\ps k_{1}\pd \ , \]
\[  \beta_{n} : \nK{n} \times_{\nG{0}} \nG{1} \rmap \nG{0} , \ \ps k_{1},\ .\ .\ .\ k_{n}, g\pd \mapsto s\ps g\pd \ . \] 
Notice that any $\alpha_n$ is \'etale. For any $n \geq 0$ we set:
\[ B_{n}(\f  ; \A ) = (\beta_{n})_{\, !}\alpha_{n}^{\, * } \A \ .\]
By \ref{losr}, the stalk at $x\in \nG{0}$ is described by:
\begin{equation}
\label{rel5.1} 
B_{n}(\f ; \A )_{x} = \gc (\beta_{n}^{-1} \ps x\pd ; \alpha_{n}^{\, *} \A ) = B_{n}( x/\f ; \pi_{x}^{\, *} \A )\ .
\end{equation}
This gives us the (stalk-wise) definition of the simplicial structure on $B_{n}(\f ; \A)$. To check the continuity, let us just remark that the boundaries can be described globally. Indeed, using the maps:
\[ d_{\, i}:  \nK{n} \times_{\nG{0}} \nG{1} = (\nG{1} /\f )^{\ps n\pd } \rmap \nK{n-1} \times_{\nG{0}} \nG{1} = (\nG{1} /\f )^{\ps n-1\pd }\]
coming from the nerve of $\nG{0} /\f $ (see \ref{nerve}, \ref{comma}), we have $\beta_{n} = \beta_{n-1} d_{\, i}$ (for all $0 \leq i \leq n$) and there are evident maps $\alpha_{n}^{\, *} \A \rmap d_{\, i}^{\, *} \alpha_{n-1}^{\, *} \A $ (compare to the definition of (\ref{rel4.2})); the boundaries of $B_{\cdot}(\f ; \A )$ are in fact: 
\[ (\beta_{n})_{\, !} \alpha_{n}^{\, *} \A = (\beta_{n-1})_{\, !} (d_{\, i})_{\, !} \alpha_{n}^{\, *} \A \rmap (\beta_{n-1})_{\, !} 
(d_{\, i})_{\, !} d_{\, i}^{\, *} \alpha_{n-1}^{\, *} \A \rmap (\beta_{n-1})_{\, !} \alpha_{n-1}^{\, *} \A  \ .\]
To describe the action of $\G$ on $B_{\cdot}(\f ; \A )$, let $g: x\rmap x\, ' $ be an arrow in $\G$. The homomorphism (\ref{rel5.0}) induces an obvious map $B_{\cdot}(x\, '/\f ; \pi_{x'}^{*}\A ) \rmap B_{\cdot}(x/\f ; \pi_{x}^{*}\A )$ which, via (\ref{rel5.1}), is the action by $g\ ( :B_{\cdot}(\f ; \A )_{x'} \rmap B_{\cdot}(\f ; \A )_{x}).$
}
\end{num}

\hspace*{.1in} If $\es$ is a c-soft $\Ka$-sheaf, $\el \f_{\, !}\es$ is defined as the chain complex of $\G$-sheaves (associated to the simplicial complex) $B_{\cdot}(\f; \es)$. If $\es_{\cdot}$ is a bounded below chain complex of c-soft $\Ka$-sheaves, define $\el \f_{\, !}\es_{\cdot}$ as the total complex of $B_{\cdot}(\f; \es_{\cdot})$. For an arbitrary $\Ka$-sheaf $\A$, $\el \f_{\, !}\A$ is defined to be $B_{\cdot}(\f; \es^{-\cdot})$ where $\es^{\cdot}$ is a resolution of $\A$ as in \ref{defhom}.  More generally, we define $\el \f_{\, !}\A_{\cdot}$ for any bounded below chain complex of $\Ka$-sheaves  using a resolution $\A_{\cdot} \rmap \R_{\cdot}$ as in \ref{hyperhom}. As in the case of homology (cf. \ref{defhom}) , we see that $\el \f_{\, !}$ is well defined up to quasi-isomorphism; in  particular, the ``derived functors'':
\[ L_{n}\f_{\, !}(-) = \Ha_{n}(\el \f_{\, !}(-)) : \abK \rmap \abG \]
are well defined up to isomorphism.  For $n= 0$ we simply denote $L_{0}\f_{\, !}: \abK \rmap \abG$ by $\f_{\, !}$. 

\begin{prop}\label{stalks} For any $x \in \nG{0}$, there are isomorphisms:
\begin{equation}
\label{rel5.4}
 (L_{n}\f_{\, !}(\A))_{x} \cong H_{n}( x/\f; \pi_{x}^{*}\A ) \ \ \ \ {\rm for \ all} \ \ x \ \in \nG{0} \ \ \ .
\end{equation}
\end{prop}

\emph{Proof:} This is an immediate consequence of relation (\ref{rel5.1}), and the fact that $\alpha_{n}^{\, *}$'s preserve c-softness since they are induced by \'etale maps.  $\Boxe$

\begin{st}
\label{spseq}
(``Leray-Hochschild-Serre spectral sequence'') For any homomorphism $\f: \Ka \rmap \G$ between \'etale groupoids and any $\Ka$-sheaf $\A$ there is a natural spectral sequence:
\[ E_{p,q}^{2} = H_{p}(\G; L_{q}\f_{\, !}\A ) \Longrightarrow H_{p+q}(\Ka; \A)\ .\]
\end{st}

\emph{Proof:} The spectral sequence follows from an isomorphism:
\begin{equation}
\label{rel5.2}
\mathbb{H}_{\, *}(\G; \el\f_{\, !}\A) \cong H_{*}(\Ka; \A) 
\end{equation}
and \ref{prop1} applied to $\el\f_{\, !}(\A)$.\\
\hspace*{.3in} To prove (\ref{rel5.2}) we consider the double complex $\C_{p,q}(\A) = B_{p}(\G; B_{q}(\f; \A))$ and we show that there are maps $\C_{0,q}(\A) \rmap B_{q}(\Ka; \A)$, functorial in $\A$, such that the augmented complex
\begin{equation}
\label{rel5.3}
 .\ .\ .\ \rmap \C_{2,q}(\A) \rmap \C_{1,q}(\A) \rmap \C_{0,q}(\A) \rmap B_{q}(\Ka; \A) \rmap 0 
\end{equation}
is acyclic for any c-soft $\Ka$-sheaf $\A$.\\
\hspace*{.3in} Using the diagram:
\[ \xymatrix{
\nK{q}\times_{\nG{0}} \nG{p+1} \ar[d]_-{v}\ar[r]^-{u} & \nG{p} \ar[d]^{\tau_{p}}\\
\nK{q}\times_{\nG{0}} \nG{1} \ar[d]_-{w}\ar[r]^-{\beta_{q}}\ar[rd]^{\alpha_{q}} & \nG{0} \\
\nK{q} \ar[r]^-{\tau_{q}} & \nK{0} } \]
where $\alpha_{q}, \beta_{q}, \tau_{q}, \tau_{p}$ are those defined before, $v, w$ are the projections into the first components, $u$ is the projection into the last components and $\gamma_{q}= wv$, we have by the general properties of the Appendix:
\begin{eqnarray}
\C_{p,q}(\A)    & = & \gc (\nG{p} ; \tau_{p}^{*} (\beta_{q})_{\, !} \alpha_{q}^{*} \A ) \nonumber\\
                & = & \gc (\nG{p} ; u_{\, !} v^{*} \alpha_{q}^{*} \A ) \nonumber\\
                & = & \gc (\nK{q} \times_{\nG{0}} \nG{p+1} ; v^{*} \alpha_{q}^{*} \A ) \nonumber\\
                & = & \gc (\nK{q} ; \ps \gamma_{q}\pd_{\, !} v^{*} \alpha_{q}^{*} \A ) \nonumber\\
                & = & \gc (\nK{q} ; \ps \gamma_{q}\pd_{\, !} \ps \gamma_{q}\pd^{*} \tau_{q}^{*} \A )\ ,\nonumber\\ 
B_{q}(\Ka; \A)  & = & \gc(\nK{q}; \tau_{q}^{*}\A) \ . \nonumber
\end{eqnarray}
\hspace*{.3in} Via these equalities, the augmented chain complex (\ref{rel5.3}) commies from an augmented simplicial sheaf on $\nK{q}$ whose stalk at
 $x \stackrel{k_{1}}{\lmap}\ .\ .\ .\ \stackrel{k_{q}}{\lmap} y$ has the form:
\[    .\ .\ .\ \rmap \bigoplus_{\f\ps y\pd \stackrel{f}{\lmap} x_{0} \stackrel{g_{1}}{\lmap} x_{1} \stackrel{g_{2}}{\lmap} x_{2}}  \A_{x}   \rmap  \bigoplus_{\f\ps y\pd \stackrel{f}{\lmap} x_{0} \stackrel{g_{1}}{\lmap} x_{1}}   \A_{x} \rmap \bigoplus_{\f\ps y\pd \stackrel{f}{\lmap} x_{0}}  \A_{x} \rmap \A_{x}\ .  \]
\hspace*{.3in} This is in fact the augmented bar complex computing the homology of the (contractible, discrete) category $\G/\f\ps y\pd$ with constant coefficients $\A_{x}$. In particular it is acyclic (with the usual contraction $(f, g_{1},\ . . . , g_{n}; a) \mapsto (1, f, g_{1},\ . . . , g_{n}; a)$).\ \ $\Boxe$

\begin{num}
\label{wequiv}
\emph{
{\bf Remarks and examples.} \\
\hspace*{.3in} 1). The isomorphism (\ref{rel5.2}) is actually a consequence of the quasi-isomorphism $\el \f_{\, !} pt_{\, !}= pt_{\, !}$ (where $pt$ is the map into the trivial groupoid); this is a particular case of the naturality property $\el\f_{\, !} \el\psi_{\, !}= \el(\f\compose\psi)_{\, !}$ (``up to quasi-isomorphism''), which can be proved in an analogous way. Compare to \cite{We1}. \\
\hspace*{.3in} 2). If $\f: \Ka \rmap \G$ is  \'etale , $\es \in \abK$, then there is no need of c-soft resolutions to define $\el\f_{\, !}\es$. Indeed, the condition on $\f$ implies that the maps $\beta_n$ defined in \ref{flosr} are \'etale, so there is a quasi-isomorphism $\el\f_{\, !}\es \simeq B_{\cdot}(\f; \es)$.\\
\hspace*{.3in} 3). Let $\f: \Ha \rmap \G$ be a morphism for which all the squares in (\ref{rel4.6}) are pullbacks. Recall that in this case, the functor $(\f_{0})_{\, !}: \underline{Ab}(\nK{0}) \rmap \underline{Ab}(\nG{0})$ ``extends'' to a functor $\f_{\, !}: \abK \rmap \abG$, making the diagram:
}
\end{num}
\[ \xymatrix {
\abK  \ar[r]^-{{\rm forget}} \ar[d]_-{\f_{\, !}} & \underline{Ab}(\nK{0}) \ar[d]^-{(\f_{0})_{\, !}} \\
\abG  \ar[r]^-{{\rm forget}} & \underline{Ab}(\nG{0}) } \]
commute. This simple minded functor of \ref{functor} agrees (up to quasi-isomorphism) with the functor $\el\f_{\, !}$, described in \ref{flosr}. Indeed, for such a morphism $\f$ and a point $x\in \nG{0}$ the comma groupoid $x/\f$ is a space (or more precisely, equivalent to the groupoid corresponding to a space, cf. \ref{ex1}.1). In this case, the spectral sequence \ref{spseq}  degenerates for c-soft sheaves $\B$ (but not for arbitrary sheaves). If $\f$ is moreover \'etale, it does always degenerate, and yields the isomorphism already proved in \ref{functor}.3.\\

\begin{cor}
\label{morinv}
(``Morita invariance'') For any Morita equivalence $\f : \Ka \rmap \G$ and any $\G$-sheaf $\A$ there is a natural isomorphism
\[ H_{p}(\G; \A) \cong H_{p}(\Ka; \f^{*}\A) .\]
\end{cor}

\emph{Proof}: Theorem \ref{spseq} gives a spectral sequence $H_{p}(\G; L_{q}\f_{\, !}\f^{*}\A) \Longrightarrow H_{p+q}(\Ka; \f^{*}\A)$. By (\ref{rel5.4}) the stalk of $L_{q}\f_{\, !}\f^{*}\A$ at a point $x\in \nG{0}$ computes the homology of the nerve of $x/\f$. If $\f$ is a Morita equivalence, this nerve is a contractible simplicial set. Thus, the spectral sequence degenerates to give an isomorphism:
\[ H_{p}(\G; L_{0}\f_{\, !}\f^{*}\A) \cong H_{p}(\Ka; \f^{*}\A)\ .\]
It thus suffices to observe that the $\G$-sheaf $L_{0}\f_{\, !}\f^{*}\A$ is isomorphic to $\A$ itself.\ \ $\Boxe$

\begin{num}
\label{fibprod}
\emph{
{\bf Fibered products of groupoids.} For homomorphisms $\f:\Ha \rmap \G$ and $\psi:\Ka \rmap \G$, their fibered product $\Ha\times_{\G}\Ka$: 
\[ \xymatrix{
\Ha\times_{\G}\Ka \ar[r]^-{q}\ar[d]_-{p} & \Ka \ar[d]^-{\psi}\\
\Ha \ar[r]^-{\f} & \G } \]
is constructed as follows. The space of objects is the space $\nH{0}\times_{\nG{0}}\nG{1}\times_{\nG{0}}\nK{0}$ of triples $\ps y, g, z\pd$ where $y\in \nH{0}, z\in \nK{0}$ and $g: \f\ps y\pd \rmap \psi\ps z\pd$ in $\G$. An arrow $\ps y, g, z\pd \rmap \ps y\, ', g\, ', z\, '\pd$ is a pair of arrows $h: y \rmap y\, '$ in $\Ha$ and $k: z\rmap z\, '$ in $\Ka$ such that $g\, '\compose \f\ps h\pd = \psi\ps k\pd \compose g$. The groupoid $\Ha\times_{\G}\Ka$ is again \'etale if $\G, \Ha, \Ka$ are. This notion of fibered product is the appropriate one for groupoids and (generalized) morphisms described in \ref{bundles} and \ref{morita}. In particular, if $\psi: \Ka \rmap \G$ is a Morita equivalence, then so is $p: \Ha\times_{\G}\Ka \rmap \Ha$.
}
\end{num}

\begin{prop} (Change-of-base formula) Consider a fibered product of \'etale groupoids as in \ref{fibprod}. For any (c-soft) $\Ka$-sheaf $\es$, there is a canonical quasi-isomorphism:
\[ \f^{*}\  \el \f_{\, !} (\es) \simeq \el p_{\, !} q^{*} (\es)\ .\]
\end{prop}

\emph{Proof:} For $y\in \nH{0}$, the comma groupoid $y_{0}/p$ is Morita equivalent to the comma groupoid $\f\ps y_{0}\pd/\psi$, (by a Morita equivalence $y_{0}/p \rmap \f\ps y_{0}\pd/\psi$ which is continuous in $y_{0}$ and which respects the action by $\Ha$). Using this observation, the proposition follows in a straightforward way from \ref{morinv} and \ref{pull}.\ \ $\Boxe$

\begin{num}
\label{comcoh}
\emph{
{\bf Compactly supported cohomology.} It is sometimes more convenient to re-index the homology groups and to see them as compactly supported cohomology groups. Because of this, we define:
\[ H_c^{n}(\G; - ) = H_{-n}(\G; - )\ \ \]
(which give a precise meaning to ``$H_c^{*}(B\G; \A)$'').The same applies to the functors $L_n\f_{\, !}$ introduced in this paragraph: if $\f: \Ka \rmap \G$ is a homomorphism, we define $R^n\f_{\, !}:= L_{-n}\f_{\, !}: \abK \rmap \abG$. With these notations, Leray spectral sequence becomes a (cohomological) spectral sequence with $E^2$-term $H_c^p(\G; R^{q}\f_{\, !}\A) \Longrightarrow H_c^{p+ q}(\Ka; \A)$. If the fibers $x/\f$ are oriented $k$-dimensional manifolds, the transgression of this spectral sequence will give ``the integration along the fibers'' map:
\[   \int_{\rm fiber} : H_c^*(\G; \mathbb{R}) \rmap H_c^{*- k}(\Ka; \mathbb{R}) .\]
}
\end{num}

\begin{num}
\emph{
\label{orbif}
{\bf Orbifolds.} As we have already mentioned in \ref{ex1}.5, orbifolds are characterized by \'etale groupoids which are proper. Let $\G$ be such a groupoid. The ``leaf space'' $M$ of $\G$ (i.e. the space obtained from $\nG{0}$ dividing out by the equivalence relation $x \sim y$ iff there is an arrow in $\G$ from $x$ into $y$), will be a Hausdorff space; it is the underlying space of the orbifold induced by $\G$ (see \cite{MP}). The obvious projection $\pi :\G\rmap M$ induces a spectral sequence:
\[ H_p(M; L_q\pi_{\, !}\A) \Longrightarrow H_{p+ q}(\G;\A), \]
for any $\A\in \abG$. The stalk of $L_q\pi_{\, !}$ at $x\in M$ is:
\[ (L_q\pi_{\, !})_{x} \cong H_q(G_{\tilde{x}}; \A_{\tilde{x}}),\]
where $\tilde{x}\in\nG{0}$ is any lift of $x$, and $G_{\tilde{x}}$ is the (finite) group $\{ \gamma\in\nG{1}: s\ps\gamma\pd = t\ps\gamma\pd =\tilde{x}\}$ (this follows from \ref{stalks} and the Morita equivalence $x/\pi \sim G_{\tilde{x}}$.\\
\hspace*{.3in} In particular, for $\A\in \underline{Mod}_{\mathbb{R}}(\G)$, the spectral sequence degenerates and gives an isomorphism:
\begin{equation}
\label{s1}
H_{*}(\G; \A) \cong H_c^{-*}(M; \pi_{\, !}\A) .
\end{equation}
\hspace*{.3in} This also shows that the ``co-invariants functor'':
\[ \Gamma_{\G}( -):= H_0(\G; -) : \underline{Mod}_{\mathbb{R}}(\G) \rmap \underline{Mod}_{\mathbb{R}} \]
is left exact and that $H_c^*(\G; -)$ (see \ref{comcoh}) are the right derived functors of $\Gamma_{\G}$.
}
\end{num}

\begin{num}
\emph{
{\bf Basic cohomology.} Let $\G$ be a smooth \'etale groupoid. The space $\Omega^*_{c, basic}(\G)$ of compactly supported basic forms is defined as the Cokernel of:
\[ \Omega^*_c(\nG{1}) \stackrel{d_0-d_1}{\rmap} \Omega^*_c(\nG{0}), \]
where $d_0, d_1$ are the maps coming from the nerve of $\G$. In other words, $\Omega^*_{c, basic}(\G)= \Gamma_{\G}(\Omega^*)$. The {\it basic compactly supported cohomology} of $\G$, denoted $H_{c, basic}^*(\G)$, is defined as the cohomology of the complex $\Omega^*_{c, basic}(\G)$ (with the differential induced by DeRham differential on $\Omega_c^*(\G)$. There is an obvious projection from the reindexed homology (see \ref{comcoh}):
\[ H_c^*(\G; \mathbb{R}) \rmap H_{c, basic}^*(\G), \]
which is an isomorphism if $\G$ is proper (cf. \ref{orbif}). In this case we also have:
\[ \Omega_{c, basic}^*(\G) \cong \{ \omega\in\Omega^*(\nG{0}): \omega \ {\rm is}\ \G-{\rm invariant},{\rm and} \ \pi(supp \, \omega) \ {\rm is\ compact\ in}\ M \ ]\]
(where $\pi: \G\rmap M$ is the projection considered in \ref{orbif}). This map associates to $\omega\in \Omega^*_{c, basic}(\G)$ the $\G$-invariant form $\tilde{\omega}$ on $\nG{0}$, given by:
\[ \tilde{\omega}\ps x\pd = \sum_{x \stackrel{g}{\rmap} y} \omega\ps y\pd g\ .\]
}
\end{num}

\section{Verdier duality}

\hspace*{.3in} In this section all sheaves are sheaves of $\mathbb{R}$- modules, i.e. real vector spaces (we can actually use any field of characteristic $0$), and $Hom$ and $\otimes$ are all over $\mathbb{R}$. We will establish a Verdier type duality for the functor $\el \f_{\, !}$ (i.e. $\f_{\, !}$ viewed at the level of the derived categories) and an associated functor $\f^{!}$ to be described, by extending one of the standard treatments \cite{Iv} to \'etale groupoids. (But our presentation is self-contained.) As a special case, we will obtain a Poincar\'e duality between the (Haefliger) cohomology of \'etale groupoids described in Section $2$ and the homology theory (Section $4$).

\begin{num}
\label{tens}
\emph{
{\bf Tensor products.} As a preliminary remark, we observe the following properties of tensor products over $\mathbb{R}$. First, if $\A$ is a c-soft sheaf on a space $Y$ and $\B$ is any other sheaf, the tensor product $\A\otimes\B$ is again c-soft. Moreover, for the constant sheaf associated to a vector space $V$ we have $\gc (Y; \A\otimes V) = \gc (Y; \A) \otimes V$ (cf. \ref{apt}). It follows by comparing the stalks that for a map $f: Y\rmap X$, also:
\[ f_{\, !}(\A \otimes \f^{*}\B) = f_{\, !}(\A) \otimes \B \]
for any sheaf $\B$ on $X$ (see \cite{Borel, Iv}). These properties extend to a morphism $\f: \Ka\rmap \G$ of \'etale groupoids: for a c-soft $\Ka$-sheaf $\A$ and any $\G$-sheaf $\B$, there is an isomorphism:
\[ \f_{\, !}(\A \otimes \f^{*}\B) = \f_{\, !}(\A) \otimes \B \ .\]
}
\end{num}

\begin{num}
\emph{
{\bf The sheaves $\mathbb{R}[V]$.} Let us fix an \'etale groupoid $\Ka$. Any open set $V \subset \nK{0}$ gives a free $\Ka$-sheaf (see \ref{ex2}.3) of sets $\tilde{V}$, given by the \'etale map $s: t^{-1}(V) \rmap \nK{0}$ and the $\Ka$-action defined by composition. Let $\mathbb{R}[V]$ be the free $\mathbb{R}$-module on this $\Ka$-sheaf $\tilde{V}$. So $\mathbb{R}[V]$ is a $\Ka$-sheaf of vector spaces, and for any other such $\Ka$-sheaf $\B$ we have:
\begin{equation}
\label{rel6.1}
Hom_{\Ka, \mathbb{R}}(\mathbb{R}[V], \B) = Hom_{\Ka}(\tilde{V}; \B) = Hom_{\nK{0}}(V, \B) = \Gamma (V; \B)
\end{equation}
(These four occurrences of $\B$ denote $\B$ as a $\Ka$-sheaf of vector spaces, as a $\Ka$-sheaf of sets, and (twice) as a sheaf on $\nK{0}$, respectively.)\\
\hspace*{.3in} There is a natural morphism:
\begin{equation}
\label{rel6.2}
e= e_{V} : \Ka/V \rmap \Ka
\end{equation}
of \'etale groupoids (of the kind described in \ref{functor}.3), and $\mathbb{R}[V]$ can also be obtained from the constant sheaf $\mathbb{R}$ on $\Ka/V$ as:
\begin{equation}
\label{rel6.3}
\mathbb{R}[V] = e_{\, !}( \mathbb{R})\ .
\end{equation}
From this point of view, the mapping properties (\ref{rel6.1}) follow by the adjunction between $e_{\, !}$ and $e^{*}$, together with the Morita equivalence $\Ka/V \simeq V$ (where $V$ is viewed as a trivial groupoid, \ref{ex1}.3).\\
\hspace*{.3in} If $V, W \subset \nK{0}$ are open sets and $\sigma: V \rmap \nK{1}$ is a section of $s: \nK{1} \rmap \nK{0}$ such that $t\compose \sigma (V) \subset W$, then composition with $\sigma$ gives a morphism $\mathbb{R}[V] \rmap \mathbb{R}[W]$. In this sense, the construction is functorial in $V$.
}
\end{num}

\begin{lem}
\label{lem6.3}
For any $\Ka$-sheaf of vector spaces $\A$ there is an exact sequence of the form:
\[ \bigoplus_{j} \mathbb{R}\, [V_{j}] \rmap \bigoplus_{i} \mathbb{R}\, [V_{i}] \rmap \A \rmap 0 \ .\]
\end{lem}

\emph{Proof:} It suffices to prove that any $\Ka$-sheaf can be covered by $\Ka$-sheaves of the form $\mathbb{R}\, [V]$, and this is clear from (\ref{rel6.1}).\ \ $\Boxe$

\begin{num}
\emph{
{\bf The sheaves $\es_{V}$.} Let $\es$ be any c-soft $\Ka$-sheaf. We write $\es_{V}$ for the sheaf $\es \otimes \mathbb{R}\, [V]$. Note that: 
\[ \es_{V} = \es \otimes e_{\, !}(\mathbb{R}) = e_{\, !}(e^{*}(\es) \otimes \mathbb{R}) = e_{\, !} e^{*}(\es)\]
(see \ref{tens}). In particular, $\es_{V}$ is again c-soft, and has the following mapping properties:
\begin{equation}
\label{rel6.4}
Hom_{\Ka}( \es_{V}, \A)\, =\, Hom_{\Ka}( \mathbb{R}\, [V],\ \underline{Hom}(\es, \A)) = \Gamma (V, \ \underline{Hom}(\es, \A)) = Hom_{V}(\es|_{V}, \A|_{V}).
\end{equation}\\
\hspace*{.3in} Now suppose $V= \bigcup V_{i}$ is an open cover. We claim that the associated sequence 
\begin{equation}
\label{rel6.5}
.\ .\ .\ \rmap \bigoplus \es_{V_{i_{0}i_{1}}}\ \  \rmap \bigoplus \es_{V_{i_{0}}}\ \ \rmap \es_{V} \rmap 0
\end{equation}
is exact. To see this, it suffices to prove that the sequence obtained by homming into any injective $\Ka$-sheaf $\ii$,
\[ 0 \rmap Hom_{\Ka}(\es_{V}, \ii) \rmap Hom_{\Ka}(\bigoplus \es_{V_{i_{0}}}, \ii) \rmap \ .\ .\ . \]
is exact. This is clear from the mapping properties (\ref{rel6.4}).
}
\end{num}

\begin{num}
\emph{
{\bf The sheaves $\f_{\, !}(\es_{V})$.} From now on let $\f :\Ka \rmap \G$ be a homomorphism between \'etale groupoids. For an open set $V \subset \nK{0}$, $\f$ induces a map $\f_{\, V}: V\rmap \G$, which fits into a commutative diagram:
}
\end{num}
\[ \xymatrix{
V \ar[d]_-{\f_{\, V}} \ar[r]^-{\sim}_-{i} & \Ka/V \ar[d]^-{e}\\
\G & \Ka \ar[l]_-{\f} . } \]
(where $i$ is the canonical Morita equivalence). Thus, for any c-soft $\Ka$-sheaf $\es$, we have:
\begin{equation}
\label{half}
 \f_{\, !}(\es_{V}) = \f_{\, !} e_{\, !} e^{*}(\es) = (\f_{\, V})_{\, !}(\es|_{V})\ .
\end{equation}
\hspace*{.3in} Notice that the groupoid $x/\f_{\, V}$ is a space (\ref{ex1}.3) for any object $x\in \nG{0}$; this and the general description of $\f_{\, !}$ (see (\ref{rel5.4})) give a simple description of the stalks of $\f_{\, !}(\es_{V})$. It follows from this description and the corresponding fact for spaces that $\f_{\, !}$ maps the exact sequence (\ref{rel6.5}) into an exact sequence:
\begin{equation}
\label{rel6.6}
.\ .\ .\ \rmap \bigoplus \f_{\, !}(\es_{V_{i_{0}i_{1}}})\ \  \rmap \bigoplus \f_{\, !}(\es_{V_{i_{0}}})\           \rmap \f_{\, !}(\es_{V}) \rmap 0\ .
\end{equation}

\begin{num}
\emph{
{\bf The $\Ka$-sheaves $\f^{\, !}(\es, \ii)$.}  Again, let $\f : \Ka \rmap \G$ be any homomorphism between \'etale groupoids, let $\es$ be a c-soft $\Ka$-sheaf, and let $\ii$ be an injective $\G$-sheaf. Define for each open set $V \subset \nK{0}$:
\[ \f^{\, !}(\es, \ii)(V) = Hom_{\G}(\f_{\, !}(\es_{V}), \ii)\ .\]
\hspace*{.3in} We claim that this defines a sheaf $\f^{\, !}(\es, \ii)$ on $\nK{0}$. Indeed, for an inclusion $V\subset W$ there is an evident map $\f^{\, !}(\es, \ii)(W) \rmap \f^{\, !}(\es, \ii)(V)$ induced by the map $\es_{V} \rmap \es_{W}$. And for a covering $V= \bigcup V_{i}$, the sheaf property follows from the injectivity of $\ii$ together with the exact sequence (\ref{rel6.6}). Furthermore, this sheaf $\f^{\, !}(\es, \ii)$ carries a natural $\Ka$-action: for any arrow $k: y\rmap z$ in $\Ka$, let $W_{y}$ and $W_{z}$ be neighborhoods of $y$ and $z$ so small that $s: \nK{1} \rmap \nK{0}$ has a section $\sigma$ through $k$ with $t\compose \sigma : W_{y}\rmap W_{z}$. Then $\sigma$ gives a map $\mathbb{R}\, [W_{y}] \rmap \mathbb{R}\, [W_{z}]$  and hence $\es_{W_{y}} \rmap \es_{W_{z}}$. By composition, one obtains a map $\f^{\, !}(\es, \ii)(W_{z}) \rmap \f^{\, !}(\es, \ii)(W_{y})$, and hence by taking germs an action
$(-)\cdot k : \f^{\, !}(\es, \ii)_{z} \rmap \f^{\, !}(\es, \ii)_{y}\ .$
}
\end{num}

\begin{prop}
\label{duality}
(Duality formula) Let $\f : \Ka\rmap \G$ be a morphism of \'etale groupoids. For any injective $\G$-sheaf $\ii$, any c-soft $\Ka$-sheaf $\es$ and any other $\Ka$-sheaf $\A$, there is a natural isomorphism of abelian groups:
\[ Hom_{\Ka}(\A, \f^{\, !}(\es, \ii)) \ \cong \ Hom_{\G}(\f_{\, !}(\A \otimes \es), \ii)\ .\]
In particular, $\f^{\, !}(\es, \ii)$ is again injective.
\end{prop}
\emph{Proof:} By \ref{lem6.3} and the fact that $\f_{\, !}$ is right exact on sequences of c-soft sheaves, it suffices to define a natural isomorphism:
\[ Hom_{\Ka}(\mathbb{R}\, [V], \f^{\, !}(\es, \ii)) \cong Hom_{\G}(\f_{\, !}(\mathbb{R}\, [V] \otimes \es), \ii)\ .\]
But, using (\ref{rel6.1}) and the definitions: 
\[ Hom_{\Ka}(\mathbb{R}\, [V], \f^{\, !}(\es, \ii)) \cong \Gamma (V, \f^{\, !}(\es, \ii)) \cong Hom_{\G}(\f_{\, !}(\es_{V}), \ii) \cong Hom_{\G}(\f_{\, !}(\mathbb{R}\, [V] \otimes \es), \ii).\ \Boxe\]
\newline
\hspace*{.3in} As for spaces \cite{Iv}, one can state and prove a somewhat stronger version of \ref{duality}, using the ``internal hom'' (see \ref{Inthom}):

\begin{prop}
\label{strduality}
(Duality formula, strong form) For any $\f, \A, \es$ and $\ii$ as in \ref{duality} there is a natural isomorphism of $\G$-sheaves
\[ \f_{*} \underline{Hom}_{\, \Ka}(\A, \f^{\, !}(\es, \ii)) \ \ \cong \ \ \underline{Hom}_{\, \G}(\f_{\, !}(\A \otimes \es), \ii)\ .\]
\end{prop}
\emph{Proof:} It suffices to prove that for any $\G$-sheaf $\B$ there is an isomorphism:
\[ Hom_{\G}(\B, \underline{Hom}_{\, \Ka}(\A, \f^{\, !}(\es, \ii))) \cong Hom_{\G}(\B, \underline{Hom}_{\, \G}(\f_{\, !}(\A \otimes \es), \ii))\ ,\]
natural in $\B$. This is immediate from \ref{duality} and \ref{tens}:
\begin{eqnarray*}
Hom_{\G}(\B, \f_{*}\underline{Hom}_{\, \Ka}(\A, \f^{\, !}(\es, \ii))) & = & Hom_{\Ka}(\f^{*}\B, \underline{Hom}_{\, \Ka}(\A, \f^{\, !}(\es, \ii))) \\
  & = & Hom_{\Ka}(\f^{*}(\B)\otimes \A, \f^{\, !}(\es, \ii)) \\
  & = & Hom_{\G}(\f_{\, !}(\f^{*}(\B)\otimes \A \otimes \es), \ii)\\
  & = & Hom_{\G}(\B \otimes \f_{\, !}(\A \otimes \es), \ii).\ \Boxe
\end{eqnarray*}

\begin{num}
\emph{
{\bf Remark.} Let $\f :\Ka \rmap \G$ be an \'etale morphism such that each of the squares in (\ref{rel4.6})  is a pull-back. Thus $\Ka = E\cross\G$ for some \'etale $\G$-space $E$. Then $\f_{\, !}$ has a simple description as in (\ref{functor}.3), and is left adjoint to $\f^{*}$. Thus:
\[ Hom_{\G}(\f_{\, !}(\A \otimes \es), \ii) = Hom_{\Ka}(\A \otimes \es, \f^{*}\ii) = Hom_{\Ka}(\A, \underline{Hom}_{\, \Ka}(\es, \f^{*}\ii)).\]
Since this holds for any $\A$, proposition \ref{duality} implies that for such a $\f$,
\[ \f^{\, !}(\es, \ii) = \underline{Hom}_{\, \Ka}(\es, \f^{*}\ii)) \ \ (\ \  = \underline{Hom}_{\, \G}(\f_{\, !}\es, \ii).\  \ ) \] 
}
\end{num}

\begin{num}
\label{ducompl}
\emph{
{\bf Duality for complexes.} We now extend these isomorphisms to (co-) chain complexes. It will be convenient to work with {\it chain} complexes for $\A$ and $\es$ and {\it cochain} complexes for $\ii$ in \ref{duality}, \ref{strduality}. Thus, we will use the following convention: if $\A$ is a chain complex and $\B$ is a cochain complex, $Hom(\A, \B)$ is the {\it cochain} complex defined by:
\[ Hom(\A, \B)^{n} = \prod_{p+q= n} Hom(\A_{p}, \B^{q})\ .\]
\hspace*{.3in} Recall for later use that if $\B^{\cdot}$ is injective and bounded below, then for any quasi-isomorphism of chain complexes $\A_{\cdot} \rmap \C_{\cdot}$ the map $Hom(\A_{\cdot}, \B^{\cdot}) \rmap Hom(\C_{\cdot}, \B^{\cdot})$ is again a quasi-isomorphism (by a standard "mapping cone" argument it is enough to prove the assertion for $\C_{\cdot}= 0$; in this case remark that $Hom(\A_{\cdot}, \B^{\cdot})$ is the total complex of a double cochain complex whose rows $Hom(\A_{\cdot}, \B^{p})$ are acyclic by the injectivity of $\B^{p}$) .\\
\hspace*{.3in} Similarly, for a bounded below chain complex $\es_{\cdot}$ of sheaves as in \ref{duality} we define the cochain complex $\f^{\, !}(\es_{\cdot}, \ii^{\cdot})$ by:
\[ \f^{\, !}(\es_{\cdot}, \ii^{\cdot}) = \prod_{p+q= n} \f^{\, !}(\es_{p}, \ii^{q})\ .\]
\hspace*{.3in} With these conventions, \ref{duality} gives an isomorphism of cochain complexes
\begin{equation}
\label{rel6.60}
 Hom_{\Ka}(\A_{\cdot}, \f^{\, !}(\es_{\cdot}, \ii^{\cdot})) \cong Hom_{\G}(\f_{\, !}(\A_{\cdot}\otimes \es_{\cdot}), \ii^{\cdot})\ ,
\end{equation}
for any cochain complex $\ii^{\cdot}$ of injective $\G$-sheaves, and any bounded below chain complexes $\A_{\cdot}$ and $\es_{\cdot}$ of $\Ka$-sheaves with $\es_{\cdot}$ c-soft. There is also an obvious ``strong'' version of (\ref{rel6.60}):
\[ \f_{*} \underline{Hom}_{\, \Ka}(\A_{\cdot}, \f^{\, !}(\es_{\cdot}, \ii^{\cdot}) \cong \underline{Hom}_{\, \G}(\f_{\, !}(\A_{\cdot}\otimes \es_{\cdot}), \ii^{\cdot})\ .\]
}
\end{num}

\begin{num}
\label{upper}
\emph{
{\bf The functor $\f^{\, !}(\ii^{\cdot})$.} Now let $d= cohdim(\nK{0})$, and fix a resolution:
\[ 0 \rmap \mathbb{R} \rmap \es^{0} \rmap\ .\ .\ .\ \rmap \es^{d} \rmap 0 \]
of the constant sheaf $\mathbb{R}$ by c-soft $\Ka$-sheaves (cf. \ref{cres}). For a cochain complex $\ii^{\cdot}$ of injective $\G$-sheaves define: 
\[ \f^{\, !}(\ii^{\cdot}) = \f^{\, !}( \es^{-\cdot}, \ii^{\cdot})\ .\]
\hspace*{.3in} Then $\f^{\, !}$ is adjoint to $\f_{\, !}$ in the derived category:
}
\end{num}

\begin{st}
\label{adj}
(Adjointness) Let $\f : \Ka \rmap \G$ be a morphism between \'etale groupoids. For any bounded below chain complex $\A_{\cdot}$ of c-soft $\Ka$-sheaves and any bounded below cochain complex $\ii^{\cdot}$ of injective $\G$-sheaves there is a natural quasi-isomorphism:
\[ Hom_{\Ka}(\A_{\cdot}, \f^{\, !}(\ii^{\cdot})) \simeq Hom_{\G}(\el \f_{\, !}(\A_{\cdot}), \ii^{\cdot})\ .\]
\end{st}
(there is also a `` strong'' version derived from \ref{strduality}).\\
 \\
{\it Proof:} Denote by $\F_{\cdot}$ the free resolution of the constant sheaf obtained by tensoring (\ref{rel2.free}) (for $\Ka$ instead of $\G$) by $\mathbb{R}$. Since $\A_{\cdot} \otimes \F_{\cdot} \rmap \A_{\cdot}$ is a quasi-isomorphism, using (\ref{rel6.60}) with $\A_{\cdot}\otimes\F_{\cdot}$ instead of $\A_{\cdot}$ and $\es_{\cdot}= \es^{-\cdot}$ , the fact that $\f^{\, !}(\ii^{\cdot})$ are injective (cf. \ref{duality}) and the general remark in \ref{ducompl} we get a quasi-isomorphism:
\begin{equation}
\label{rel6.7}
 Hom_{\Ka}(\A_{\cdot}, \f^{\, !}(\es_{\cdot}, \ii^{\cdot}) \simeq Hom_{\G}(\f_{\, !}(\A_{\cdot}\otimes \F_{\cdot}\otimes \es_{\cdot}), \ii^{\cdot})\ .
\end{equation}
\hspace*{.3in} Remark that for any bounded below chain complex $\B_{\cdot}$ of c-soft $\Ka$-sheaves we have a quasi-isomorphism followed by an isomorphism: $\el \f_{\, !}\B_{\cdot} \simeq \el \f_{\, !}(\B_{\cdot}\otimes\F_{\cdot}) = \f_{\, !}(\B_{\cdot}\otimes\F_{\cdot})$. This for $\B_{\cdot}= \A_{\cdot} \otimes \es_{\cdot}$ and the quasi-isomorphism $\A_{\cdot} \simeq \A_{\cdot} \otimes \es_{\cdot}$ give $\el \f_{\, !}\A_{\cdot} \simeq \f_{\, !}(\A_{\cdot}\otimes\F_{\cdot}\otimes\es_{\cdot})$. Using this, the fact that $\ii^{\cdot}$ 's are injective and (\ref{rel6.7}), the statement of the theorem follows easily.\ \ $\Boxe$\\
\newline
\hspace*{.3in} Poincar\'e duality follows in the usual way:

\begin{cor}
\label{poinc}
(Poincar\'e duality) Let $\Ka$ be an \'etale groupoid, and suppose that $\nK{0}$ is a topological manifold of dimension $d$. Let $or$ be the orientation $\Ka$-sheaf (\ref{ex2}.6). There is a natural isomorphism:
\[ H^{p+d}(\Ka; or) \cong H_{p}(\Ka; \mathbb{R})^{\vee}\ \ \ \  (p \in \mathbb{Z}).\]
\end{cor}

{\it Proof:} Let $\G= 1$ be the trivial groupoid. In \ref{adj}, let $\ii^{\cdot}$ be the complex $\mathbb{R}$ concentrated in degree $-p$, and let $\A_{\cdot}$ be the complex $\A_{i}= \es^{-i}$ (as in \ref{upper}). As $\A_{\cdot}$ is quasi-isomorphic to $\mathbb{R}$, the complex on the right of the quasi-isomorphism in \ref{adj} has:
\[ H^{0}( Hom(\f_{\, !}(\A_{\cdot}), \ii^{\cdot}) = H_{p}(\Ka; \mathbb{R})^{\vee}\ .\]
\hspace*{.3in} Now consider the left hand side of the quasi-isomorphism of \ref{adj}, $Hom(\A_{\cdot}, \f^{\, !}(\ii^{\cdot}))$ in this special case. Note first that $\f^{\, !}(\mathbb{R})$ is quasi-isomorphic to the orientation sheaf concentrated in degree $-d$:
\begin{equation}
\label{rel6.8}
\f^{\, !}(\mathbb{R}) \simeq or\, [d] \ \ .
\end{equation}
Indeed, for any open set $V \subset \nK{0}$, $\f^{\, !}(\mathbb{R})(V)= Hom(\Gamma_c(V; \es^{-\cdot}); \mathbb{R})$ (see (\ref{half})) is the complex which computes $H_c^{-\cdot}(V; \mathbb{R})^{\vee}$ , and the argument for (\ref{rel6.8}) is just like the one for spaces. Thus, for $\ii^{\cdot}= \mathbb{R}\, [p]$ and $\A_{\cdot}= \es^{-\cdot} \simeq \mathbb{R}\, [0]$,
\[ H^{0}( Hom_{\Ka}(\A_{\cdot}, \f^{\, !}(\ii^{\cdot})) = H^{0} \Gamma_{inv}(\Ka, \f^{\, !}(\ii^{\cdot})) = H^{p+d}(\Ka, or)\ .\ \Boxe \]

\section{Relation to cyclic homology}

In \cite{Cra} the homology of \'etale groupoids and the Leray spectral sequence (paragraph 5) are the main tools in dealing with cyclic homology of \'etale groupoids. We shall give here an overview of the main results in \cite{Cra} expressing the cyclic homology of (the convolution algebra of) an \'etale groupoid in terms of the homology of \'etale categories. For instance \ref{perrii} generalizes the previous results of Connes (for $\G$ a space, \cite{Co2}), Burghelea and Karoubi (for $\G$ a group), Brylinski and Nistor (for $\G$ a separated \'etale groupoid, \cite{BrNi}), Feigin and Tsygan (\cite{FeTs}), Nistor (\cite{Ni1}). In this section, all sheaves are sheaves of vector spaces, and we work in the category of smooth \'etale groupoids.

\begin{num}
{\bf Mixed complexes of sheaves.}\emph{Let $\G$ be an \'etale groupoid. By a mixed complex of $\G$-sheaves we mean a mixed complex $(\A_{\cdot}, b, B)$ in the category $\abG$. This means a family $\{ \A_{n}: n\geq 0\} $ of $\G$-sheaves equipped with maps of degree $-1$, $b: \A_{n} \rmap \A_{n-1}$ and maps of degree $+1$, $B: \A_{n} \rmap \A_{n+1}$, such that $b^2= B^2= (b+ B)^2= 0$. For the general notions and constructions concerning mixed complexes in any abelian category see \cite{Kassel}. Recall that any such mixed complex $(\A_{\cdot}, b, B)$ gives rise to a double complex $\B(\A)$ in $\abG$, hence the Hochschild and cyclic homology sheaves are defined (see also 3.6 in \cite{Cra}):}
\[  \widetilde{HH}_{*}(\A_{\cdot}), \widetilde{HC}_{*}(\A_{\cdot}) \ \in \ \abG \ . \]
\emph{The Hochschild and cyclic (hyper-) homology of the mixed complex $(\A_{\cdot}, b, B)$, denoted $HH_{*}(\G; \A)$, $HC_{*}(\G; \A)$ are defined as the hyperhomology of the complexes of $\G$-sheaves $(\A_{\cdot}, b)$ and $Tot(\B(\A))$, respectively (compare to \cite{Kar1, Kar2}). From \ref{prop1} we get two spectral sequences with $E^2$ terms:}
\[ H_{p}(\G;\widetilde{HH}_{q}(\A_{\cdot})) \Longrightarrow HH_{p+q}(\G;\A_{\cdot})\ \ {\rm and}\ \ H_{p}(\G;\widetilde{HC}_{q}(\A_{\cdot})) \Longrightarrow HC_{p+q}(\G;\A_{\cdot}). \]
\emph{The spectral sequences of the double complex $B_{\cdot}(\G;T_{\cdot})$  where $T_{\cdot}=(\A_{\cdot},b)$  or $Tot(\B(\A))$ give two spectral sequences with $E^{2}$-terms:}
\[ HH_{p}(H_{q}(\G;\A_{\cdot})) \Longrightarrow HH_{p+q}(\G;\A_{\cdot})  \ \  \mbox{and} \ \ HC_{p}(H_{q}(\G;\A_{\cdot})) \Longrightarrow HC_{p+q}(\G;\A_{\cdot}) .\]
\emph{Also from \ref{les} we get the SBI sequence relating $HH_{*}(\G; \A_{\cdot})$ and $HC_{*}(\G; \A_{\cdot})$. The periodic cyclic homology is defined (as ``usual'') as $lim_{S} HC_{*}(\G; \A_{\cdot})$.}
\end{num}

\begin{num}
\label{cii}
\emph{
{\bf Cyclic $\G$-sheaves.} It is well known (and it is a motivating example) that any cyclic object in an abelian category gives rise to a mixed complex (\cite{Kassel}). In particular, any cyclic $\G$-sheaf $\A_{\cdot}$ (i.e. a contravariant functor $\Lambda \rmap \abG$ from Connes' category $\Lambda$ (\cite{Co1}), induces a mixed complex of $\G$-sheaves. The corresponding homologies are still denoted $\widetilde{HH}_{*}(\A_{\cdot}), \widetilde{HC}_{*}(\A_{\cdot})$, $HH_{*}(\G; \A_{\cdot})$, $HC_{*}(\G; \A_{\cdot})$. \\
\hspace*{.3in} A basic example of a cyclic $\G$-sheaf for a smooth \'etale groupoid is $\cii$ defined by $\cii \ps n\pd =$the pullback of the sheaf of smooth functions on $(\nG{0})^{n+ 1}$ along the diagonal embedding $\delta_{n}: \nG{0} \rmap (\nG{0})^{n+ 1}$, with the cyclic structure described as follow. At $c \in \nG{0}$, the stalk of $\cii \ps n\pd$ is the vector space of germs $f\ps x_0, .\ .\ .\ , x_n \pd$ of smooth functions defined for $x_0, .\ .\ .\ , x_n \in \nG{0}$ around $c$, and:
}
\end{num}

\[ (d_i\, f)\ps x_0, .\ .\ .\ , x_{n- 1}\pd = \left\{ \begin{array}{ll}
                                                       f\ps x_0, .\ .\ .\ , x_i, x_i, .\ .\ .\ , x_{n- 1}\pd & \mbox{if $0 \leq i \leq n-1$} \\
                                                       f\ps x_{n- 1}, x_0, .\ .\ .\ , x_{n- 1}\pd & \mbox{if $i=n$}
                                                \end{array}
                                        \right. ,\]
\[ (t_n\, f)\ps x_0, .\ .\ .\ , x_n\pd = f\ps x_n, x_0, .\ .\ .\ , x_{n- 1}\pd .\]
Using the quasi-isomorphism $(\cii, b) \simeq (\Omega^*, 0)$ which appears in the work of Connes (see also lemma 3.5 in \cite{Cra}), we get:
\[ HH_{n}(\G; \cii) = \bigoplus_{p+ q= n} H_{p} (\G; \Omega^{q}), \]
\[ HP_{i}(\G; \cii) = \prod_{k} H_{i+ 2k}(\G), \ \ \  , i \in \{ 0, 1\}. \]

\begin{num}
\emph{
{\bf Cyclic homology of the convolution algebra.} The convolution algebra of a smooth \'etale groupoid $\G$ was used by Connes as a non-commutative model for the ``leaf space'' of $\G$. When $\G$ is Hausdorff, $C_{c}^{\infty}(\G)$ is the (locally convex) algebra of compactly supported smooth functions on $\nG{1}$, with the convolution product $(u\, v)\ps g\pd =\sum_{g_1g_2= g} u\ps g_1\pd v\ps g_2\pd $. Its (continuous) Hochschild and cyclic homology are computed by the cyclic vector space $\cii(\G)^{\natural}$ with $\cii(\G)^{\natural}= \cii_c(\G)^{\hat{\otimes}\ps n+ 1\pd} \cong \cii_c(\G^{\ps n+ 1\pd})$ (here $\hat{\otimes}$ denotes the projective tensor product and the last isomorphism is an algebraic one, see \cite{Gr}). \\
\hspace*{.3in} Using the functor $\Gamma_c$ described in the Appendix, the definition of the convolution algebra $\cii_c(\G)$ extends to the non-Hausdorff case (see $3.31$ and $3.32.2$ in \cite{Cra}). It also becomes clear that the (continuous version of the) Hochschild, cyclic and periodic cyclic homology of this algebra should be defined using the cyclic vector space $\cii(\G)^{\natural}$
with $\cii(\G)^{\natural}= \cii_c(\G^{\ps n+ 1\pd})$ (and the usual cyclic structure). In this way, the Chern-Connes-Karoubi character $Ch: K_i (\cii_c(\G)) \rmap HP_i(\cii_c(\G)),\ \ i\in \{ 0, 1\}$, extends to the non-Hausdorff case.
}
\end{num} 

\begin{num}
\emph{
{\bf The Hochschild homology of $\cii_c(\G)$.} For $\G$ as before (not necessarily Hausdorff), we introduce the groupoid of loops $\Omega(\G)= \nB{0}\cross_{\nG{0}}\G)$ (see \ref{ex1}.6), where $\nB{0}=\{ \gamma\in \nG{1}: s\ps\gamma\pd = t\ps\gamma\pd \}$ is the space of loops with the $\G$-action given by conjugation $(\gamma, g) \mapsto g^{-1}\gamma g$.\\
\hspace*{.3in} There is a (simplicial) complex $\cii_{tw}$ of c-soft $\Omega(\G)$-sheaves, which is obtained by twisting $\cii$ (see \ref{cii}) by loops. More precisely, $\cii_{tw}\ps n\pd= s^{*}\cii\ps n\pd$ (where $s: \nB{0} \rmap \nG{0}$ denotes the restriction to $\nB{0}$ of the source map), with the twisted boundaries $d_i^{tw}$ defined as follows. At $\gamma \in \nB{0}$, the stalk of $\cii_{tw}\ps n\pd$ is the vector space of germs $f\ps x_0, .\ .\ .\ , x_n\pd$ defined for $x_0, .\ .\ .\ , x_n \in \nG{0}$ around $s\ps \gamma\pd$, and:
\[ (d_i^{tw}\, f)\ps x_0, .\ .\ .\ , x_{n- 1}\pd = \left\{ \begin{array}{ll}
                                                       f\ps x_0, .\ .\ .\ , x_i, x_i, .\ .\ .\ , x_{n- 1}\pd & \mbox{if $0 \leq i \leq n-1$} \\
                                                       f\ps \gamma\ps x_{n- 1}\pd , x_0, .\ .\ .\ , x_{n- 1}\pd & \mbox{if $i=n$}
                                                \end{array}
                                        \right. ,\]
The following is a reformulation of $4.1, 3.36$ in \cite{Cra}:
}
\end{num}

\begin{prop}
\label{hoch} For any smooth \'etale groupoid $\G$, there is a natural isomorphism:
\[ HH_{*}(\cii_c(\G))= \mathbb{H}_{*}(\Omega(\G); \cii_{tw}).\]
\end{prop}

Recall that $\Lambda_{\infty}$ is the category defined in the same way as $\Lambda$, except that the cyclic relation $(t_n)^{n+ 1}= 1$ is omitted (\cite{FeTs}). Remark that $\cii_{tw}$ actually has a $\Lambda_{\infty}$ structure given by the boundaries $d_i^{tw}$ just described and the cyclic action:
\[ (t_n\, f)\ps x_0, .\ .\ .\ , x_n\pd = f\ps \gamma\ps x_n\pd , x_0, .\ .\ .\ , x_{n- 1}\pd .\]
In other words, $\cii_{tw}$ can be viewed as a contravariant functor $\Lambda_{\infty} \rmap Sh(\Omega(\G))$, or, equivalently, as a $\Lambda_{\infty}\times\Omega(\G)$- sheaf. The previous proposition becomes:
\begin{equation}
\label{relc1}
HH_{*}(\cii(\G)) = H_{*}(\Lambda_{\infty}\times\Omega(\G); \cii_{tw}) \ .
\end{equation}

\begin{num}
\emph{
{\bf The cyclic homology of $\cii_c(\G)$.} We relate now $HC_{*}(\cii_c(\G))$ to the homology of an \'etale category. Remark that:
\begin{equation}
\label{relc2}
((t_n)^{n+ 1}\, f)\ps x_0, .\ .\ .\ , x_n\pd = f\ps \gamma\ps x_n\pd , \gamma\ps x_0\pd , .\ .\ .\ , \gamma\ps x_{n- 1}\pd\pd
\end{equation}
which shows that $\cii_{tw}$ is in fact a $\Lambda_{\infty}\wedge \Omega(\G)$- sheaf, where $\Lambda_{\infty}\wedge \Omega(\G)$ is the \'etale category obtained from $\Lambda_{\infty}\times\Omega(\G)$ by imposing the relations $(t_n^{n+ 1}, id_{\gamma}) = (id_{[n]}, \gamma)$, for all $\gamma \in \nB{0}, n \geq 0$. A reformulation of $4.1, 3.36$ and $3.20$ in \cite{Cra} is:
}
\end{num}

\begin{prop} 
\label{relc3} For any smooth \'etale groupoid $\G$, there is a natural isomorphism:
\[ HC_{*}(\cii_c(\G)) = H_{*}(\Lambda_{\infty}\wedge \Omega(\G); \cii_{tw}) .\]
\end{prop}

We remark that the SBI sequence can be described (via the isomorphisms (\ref{relc1}), and the one of \ref{relc3}) as the Gysin-type long exact sequence arising from the Leray spectral sequence applied to the obvious projection map $\Lambda_{\infty} \times\Omega(\G) \rmap \Lambda_{\infty}\wedge \Omega(\G)$. Note that we have tacitly made use of the extension of homology for \'etale groupoids to \'etale categories (cf. \ref{etalecate}).

\begin{num}
\label{rkk}
\emph{
{\bf Remark.} A Morita equivalence $\G \tilde{\rmap} \Ka$ induces Morita a equivalence $\Omega(\G) \tilde{\rmap} \Omega({\Ka})$, hence the Morita invariance for homology \ref{morinv}, implies the Morita invariance of the Hochschild and cyclic homology of the smooth convolution algebras.
}
\end{num}

\begin{num}
\emph{
{\bf Localization.} Remark that $\G \subset \Omega(\G)$ as units, so we recover $HH_{*}(\G; \cii)$ as ``localization at units'':
\[ HH_{*}(\G; \cii) = HH_{*}(\cii(\G))_{[1]} .\]
The isomorphisms described in \ref{cii} will give:
\[ HH_n(\cii_c(\G))_{[1]} = \bigoplus_{p+ q= n} H_{p}(\G; \Omega^{\,q}) \]
(where $\Omega^{\,q}$ is the $\G$-sheaf of $q$-forms) and also (Theorem 4.3 in \cite{Cra}):
}
\end{num}

\begin{prop} 
\label{perrii} 
For any smooth \'etale groupoid $\G$, there is a natural isomorphism:
\[ HP_{i}(\cii_c(\G))_{[1]} = \prod_{k} H_{i+ 2k} (\G) \ \ , i \in \{0, 1\}.\]
 \end{prop} 

More generally, any $\G$-invariant subspace $\OO \subset \nB{0}$ defines a groupoid $\Omega_{\OO}(\G)= \OO\cross\nG{0} \subset \Omega(\G)$ and the localized Hochschild and cyclic homology (indicated by the subscript $\OO$). When $\OO$ is elliptic (i.e. $ord(\gamma) < \infty$, for all $\gamma \in \OO$), it is shown in \cite{Cra} (Theorem 4.4) that:
\begin{equation}
\label{peri}
HP_i(\cii_c(\G))_{\OO} = \prod_k H_{i+ 2k}(\Omega_{\OO}(\G))\ \ , \  i\in \{ 0, 1\} .
\end{equation}

\begin{num}
\emph{
{\bf The case of orbifolds.} Let $\mathcal{M}= (M, \mathcal{U})$ be an orbifold ($M$ is the underlying topological space, $\mathcal{U}$ an orbifold atlas). Due to remark \ref{rkk}, the Hochschild, cyclic and periodic cyclic homologies do not depend on the representation of the orbifold $\mathcal{M}$ by a smooth proper \'etale groupoid. We simply denote these homologies by $HH_*(\mathcal{M}), HC_*(\mathcal{M}), HP_*(\mathcal{M})$.\\
\hspace*{.3in} Note that for any representation of $\mathcal{M}$ by a proper \'etale groupoid $\G$, the loop groupoid $\Omega(\G)$ is again a proper \'etale groupoid. Denote by $\Omega(\mathcal{M})$ the underlying space of the orbifold induced by $\Omega(\G)$ (i.e. $\Omega(\mathcal{M})$ is the leaf space of $\Omega(\G)$). This space can be constructed directly by using an orbifold atlas for $\mathcal{M}$; it was introduced in this way in \cite{Kaw}; it serves there for the definition of a geometric Chern character, needed in the formulation of the index theorem for orbifold (from this point of view, the next proposition explains this choice). Alternatively, representing $\mathcal{M}$ as a quotient $N/L$, where $L$ is a compact Lie group acting on $M$, with finite stabilizers (see \cite{MP}), then:
\[ \Omega(\mathcal{M}) = \widehat{N}/L,\]
where $\widehat{N}=\{ \ps x, \gamma\pd\in M\times L : x\gamma = x\}$ is Brylinski's space, with the $L$-action $\ps x, \gamma\pd g = \ps xg, g^{-1}\gamma g\pd$. Then  (\ref{s1}) applied to $\Omega(\G)$, and (\ref{peri}), give the following result which also makes the connection with Kawasaki's definition of the Chern character for orbifolds (\cite{Kaw}):
}
\end{num}

\begin{prop} For any orbifold $\mathcal{M}$:
\[ HP_i(\mathcal{M}) = \prod_k H_c^{i+ 2k}(\Omega(\mathcal{M})) \ \ , i\in \{0, 1\}\ .\]
\end{prop}
\section{Appendix: Compact supports in non-Hausdorff spaces}

\hspace*{.3in} In this appendix we explain how the usual notions concerning compactness and sheaves on Hausdorff spaces extend to our more general context (see \ref{assumptions}). For basic definitions and facts for sheaves on Hausdorff spaces, we refer the reader to any of the standard sources \cite{God, Iv, Borel}.

\begin{num}
\label{c-s}
\emph{
{\bf c-soft sheaves.} Let $X$ be a space satisfying the general assumptions in \ref{assumptions}. An abelian sheaf $\A$ on $X$ is said to be {\it c-soft} if for any Hausdorff open $U \subset X$ its restriction $\A|_{U}$ is a c-soft sheaf on $U$ in the usual sense. By the same property for Hausdorff spaces, it follows that c-softness is a local property, i.e., a sheaf $\A$ is c-soft iff there is an open cover $X= \bigcup U_{i}$ such that each $\A|_{U}$ is a c-soft sheaf on $\A$.
}
\end{num}

\begin{num}
\label{fgam}
\emph{
{\bf The functor \gc .} Let $\A$ be a c-soft sheaf on $X$ and let $\A\,'$ be its Godement resolution (i.e. $\A\,'(U)= \Gamma(U_{{\rm discr}}; \A)$ is the set of all (not necessarily continuous) sections, for any open $U\subset X$). For any Hausdorff open set $W\subset X$, let $\gc (W, \A)$ be the usual set of compactly supported sections. If $W\subset U$, there is an evident homomorphism, ``extension by 0'' $\gc (W, \A) \rmap \gc (U, \A)\subset \Gamma(U, \A\, ')$. For any (not necessarily Hausdorff) open set $U\subset X$, we define $\gc (U, \A)$ to be the image of the map:
\[ \bigoplus_{W}\gc(W, \A) \rmap \Gamma(U, \A\,')\ , \]
where $W$ ranges over all Hausdorff open subsets $W\subset U$.\\
\hspace*{.3in} Observe that $\gc(U, \A)$ so defined is evidently functorial in $\A$, and that for any inclusion $U\subset U '$ we have an ``extension by zero'' monomorphism:
\[ \gc(U, \A) \rmap \gc(U', \A)\ .\]
\hspace*{.3in} The following lemma shows that in the definition of $\gc(U, \A)$ it is enough to let $W$ range over a Hausdorff open cover of $U$; in particular, it shows that the definition agrees with the usual one if $U$ itself is Hausdorff.
}
\end{num}

\begin{lem}
\label{cov}
Let $\A$ be a c-soft sheaf on $X$. For any open cover $U= \bigcup W_{i}$ where each $W_{i}$ is Hausdorff, the sequence:
\[ \bigoplus_{i}\gc(W_{i}, \A) \rmap \gc(U, \A) \rmap 0\]
is exact.
\end{lem}

\emph{Proof:} It suffices to show that for any Hausdorff open $W\subset U$, the map $\bigoplus_{i} \gc(W\bigcap W_{i}, \A) \rmap \gc(W, \A)$ is surjective. This is well known (see e.g. \cite{God}).\ \ $\Boxe$\\

\hspace*{.3in} This lemma is in fact a special case of the following Proposition (``Mayer-Vietoris''):

\begin{prop}
\label{MV}
Let $X= \bigcup_{i}U_{i}$ be an open cover indexed by an ordered set $I$, and let $\A$ be a c-soft sheaf on $X$. Then there is a long exact sequence:
\begin{equation}
\label{relA1}
 .\ .\ .\ \rmap \bigoplus_{i_{0}< i_{1}} \gc(U_{i_{0}i_{1}}, \A) \rmap \bigoplus_{i_{0}}\gc(U_{i_{0}}, \A) \rmap \gc(X, \A) \rmap 0
\end{equation}
Here $U_{i_{0}. . . i_{n}}= U_{i_{0}}\cap .\ .\ .\ \cap U_{i_{n}}$, as usual. (There is of course a similar exact sequence if $I$ is not ordered.)
\end{prop}
\emph{Proof:} The proposition is of course well known in the case where $X$ is a paracompact Hausdorff space. We first reduce the proof to the case where each of the $U_{i}$ is Hausdorff, as follows. Let $X= \bigcup_{j\in J}W_{j}$ be a cover by Hausdorff open sets, and consider the double complex:
\[ C_{p,q}= \bigoplus \gc(W_{j_{0} . . . j_{p}}\cap U_{i_{0} . .  i_{q}}, \A)\ ,\]
where the sum is over all $j_{0}< .\ .\ .\ < j_{p},\ i_{0}< .\ .\ .\ < i_{q}$. For a fixed $p \geq 0$, the column $C_{p, \cdot}$ is a sum of exact Mayer-Vietoris sequences for the Hausdorff open sets $W_{j_{0} . . . j_{p}}$, augmented by $C_{p,-1}= \bigoplus_{j_{0}< . . . < j_{p}}\gc(W_{j_{0} . . . j_{p}}, \A)$. Keeping the notation $U_{i_{0} . .  i_{q}}= X= W_{j_{0} . . . j_{p}}$ if $q= -1= p$, we observe that for a fixed $q\geq -1$, the row $C_{\cdot, q}$ is a sum of Mayer-Vietoris sequences for the spaces $U_{i_{0} . .  i_{q}}$ with respect to the open covers $\{ W_{j}\cap U_{i_{0} . .  i_{q}}\}$. So, if the proposition would hold for covers by Hausdorff sets, each row $C_{\cdot, q}\ (q\geq -1)$ is also exact. By a standard double complex argument it follows that the augumentation column $C_{-1, \cdot}$ is also exact, and this column is precisely the sequence in the statement of the proposition. This shows that it suffices to show the proposition in the special case where each $U_{i}$ is Hausdorff.\\
\hspace*{.3in} So assume each $U_{i}\subset X$ is Hausdorff. Observe first that exactness of the sequence (\ref{relA1}) at $\gc(X, \A)$ now follows by Lemma \ref{cov}. To show exactness elsewhere, consider for each finite subset $I_{0}\subset I$ the space $U^{I_{0}}= \bigcup_{i\in I_{0}} U_{i}$ and the subsequence:
\begin{equation}
\label{relA2}
 .\ .\ .\ \rmap \bigoplus_{i_{0}< i_{1} \ {\rm in} \ I_{0}} \gc(U_{i_{0}i_{1}}, \A) \rmap \bigoplus_{i_{0} \ {\rm in}\ I_{0}}\gc(U_{i_{0}}, \A) \rmap \gc(U^{I_{0}}, \A) \rmap 0
\end{equation}
of (\ref{relA1}). Clearly (\ref{relA1}) is the directed union of the sequences of the form (\ref{relA2}), where $I_{0}\subset I$ ranges over all finite subsets of $I$. So exactness of (\ref{relA1}) follows from exactness of each such sequence of the form (\ref{relA2}). Thus, it remains to prove the proposition in the special case of a {\it finite} cover $\{ U_{i}\}$ of $X$ by Haudorff open sets.\\
\hspace*{.3in} So assume $X= U_{1}\cup . . . \cup U_{n}$ where each $U_{i}$ is Hausdorff. For $n= 1$, there is nothing to prove. For $n= 2$, the sequence has the form 
\[ 0 \rmap \gc(U_{1}\cap U_{2}, \A) \rmap \gc(U_{1}, \A) \bigoplus \gc(U_{2}, \A) \rmap \gc(U_{1}\cap U_{2}, \A)\rmap 0\ .\]
This sequence is exact at $\gc (X, \A)$ by \ref{cov}, and evidently exact at other places. Exactness for $n= 3$ can be proved using exactness for $n= 2$. Indeed, consider the following diagram, whose upper two rows are the sequences for $n= 2, 3$ and whose third row is constructed by taking vertical cokernels, so that all columns are exact (we delete the sheaf $\A$ from the notation)(compare to pp. 187 in \cite{BoTu}):
\[ \xymatrix {
 0 \ar[d] & 0 \ar[d] & 0 \ar[d] & 0 \ar[d] & \\
  0 \ar[d]\ar[r] & \gc(U_{12}) \ar[d]\ar[r] & \gc(U_{1})\oplus\gc(U_{2}) \ar[d]\ar[r] &\gc(U_{1}\cup U_{2})\ar[r]\ar[d] & 0 \\
  \gc(U_{123}) \ar[d]\ar[r] & \oplus_{1\leq i< j\leq 3}\gc(U_{ij}) \ar[d]\ar[r]& \gc(U_{1}) \oplus \gc(U_{2}) \oplus \gc(U_{3}) \ar[d]\ar[r]  & \gc(U_{1}\cup U_{2}\cup U_{3}) \ar[d]\ar[r] & 0\\
  \gc(U_{123}) \ar[r]\ar[d] & \gc(U_{13})\oplus\gc(U_{23}) \ar[d]\ar[r] & \gc(U_{3}) \ar[d]\ar[r]^-{\pi} & C \ar[d]\ar[r] & 0 \\
   0 & 0 & 0 & 0 & } \]
To show that the middle row is exact, it thus suffices to prove that the lower row is exact. This row can be decomposed into a Mayer-Vietoris sequence for the case $n= 2$, already shown to be exact,
\[ 0 \rmap \gc(U_{123}) \rmap \gc(U_{13}) \oplus \gc(U_{23}) \rmap \gc(U_{3} \cap (U_{1}\cup U_{2})) \rmap 0 \]
and the sequence:
\[ 0 \rmap \gc(U_{3}\cap (U_{1}\cup U_{2})) \rmap \gc(U_{3}) \rmap C \rmap 0\ .\]
\hspace*{.3in} The exactness of the latter sequence is easily proved by a diagram chase, using exactness of the right-hand column.\\
\hspace*{.3in} An identical argument will show that the exactness for a cover by $n+1$ opens follows from exactness for one by $n$ opens, so the proof is completed by induction.\ \ $\Boxe$\\
\newline
\hspace*{.3in} Proposition \ref{MV} is our main tool for transfering standard facts from sheaf theory on Hausdorff spaces to the non-Hausdorff case, as illustrated by the following corollaries.

\begin{cor}
\label{pair}
 Let $Y\subset X$ be a closed subspace, and let $\A$ be a c-soft sheaf on $X$. There is an exact sequence
\[ 0 \rmap \gc (X-Y, \A) \stackrel{i}{\rmap} \gc(X, \A) \stackrel{r}{\rmap} \gc(Y, \A) \rmap 0\]
($i$ is extension by zero, $r$ is the restriction).
\end{cor}

\emph{Proof:} This (including the fact that the map $r$ is well defined) follows by elementary homological algebra from the fact that the Corollary holds for Hausdorff spaces, by using \ref{MV} for a cover of $X$ by Hausdorff open sets $U_{i}$ and for the induced covers of $Y$ by $\{ U_{i}\cap Y\} $ and $X-Y$ by $\{ U_{i}-Y\}$.\ \ $\Boxe$

\begin{cor}
\label{apt}
 For a family $\A_{i}$ of c-soft sheaves on $X$ the direct sum $\oplus\A_{i}$ is again c-soft, and:
\[ \gc(X, \oplus\A_{i}) \cong \oplus\gc(X, \A_{i})\ .\]
\hspace*{.3in} In particular, when working over $\mathbb{R}$, we have for any c-soft sheaf $\es$ of $\mathbb{R}$-vector spaces and any vector space $V$ that the tensor product $\es\otimes_{\mathbb{R}}V$ (here $V$ is the constant sheaf) is again c-soft, and the familiar formula:
\begin{equation}
\label{relA3}
\gc(X, \es\otimes_{\mathbb{R}}V) \cong \gc(X, \es) \otimes_{\mathbb{R}}V\ .
\end{equation}
\end{cor}

\begin{cor}
\label{qi}
Let $\A_{\cdot} \rmap \B_{\cdot}$ be a quasi-isomorphism between chain complexes of c-soft sheaves on $X$. Then:
\[ \gc(X, \A_{\cdot}) \rmap \gc(X, \B_{\cdot})\]
is again a quasi-isomorphism.
\end{cor}

\emph{Proof:} By a ``mapping cone argument'' we may assume that $\B_{\cdot}= 0$. In other words, we have to show that $\gc(X, \A_{\cdot})$ is acyclic whenever $\A_{\cdot}$ is. This follows from the Mayer-Vietoris sequence \ref{MV} together with the Hausdorff case.\\
\hspace*{.3in}(We remark that it is necessary to assume that the chain complexes are bounded below if $X$ does not have locally finite cohomological dimension, as in \ref{assumptions}).\ \ $\Boxe$\\
 \\
\hspace*{.3in} The following Corollary is included for application in \cite{Cra}.

\begin{cor} Let $Y\subset X$ be a closed subspace, and let $\theta: X\rmap \mathbb{R}$ be a continuous map such that $\theta^{-1}(0)= Y$. Let $\A$ be a c-soft sheaf on $X$. Then for any $\alpha \in \gc(X, \A)$,
\[ \alpha|_{\, Y}= 0 \ \ {\rm iff} \ \ \exists\ \varepsilon> 0:\ \ \alpha|_{\theta^{-1}(-\varepsilon, \varepsilon)} = 0 \]
(here $\alpha|_{Y}$ is the restriction $r(\alpha)$ as in \ref{pair}).
\end{cor}

\emph{Proof:} For $\varepsilon \geq 0$, write $Y_{\varepsilon}=\{ x\in X: |\theta\ps x\pd | \leq \varepsilon\} $, and for each open $U\subset X$ write
\[ \gc^{\varepsilon}(U, \A) = \{ \alpha \in \gc(U, \A): \alpha |_{U\cap Y_{\varepsilon}}= 0\}\ .\]
\hspace*{.3in} It suffices to show that: 
\[ \bigoplus_{\varepsilon > 0} \gc^{\varepsilon}(X, \A) \rmap \gc^{0}(X, \A) \]
is epi. Let $\{ U_{i}\} $ be a cover of $X$ by Hausdorff open sets, and consider the diagram:
\[ \xymatrix{
\bigoplus_{i, \varepsilon> 0} \gc^{\varepsilon}(U_{i}, \A) \ar[d]\ar[r]^-{u} & \bigoplus_{i}\gc^{0}(U_{i}, \A) \ar[d]^-{\pi}\ar[r]^-{\sim} & \bigoplus_{i}\gc(U_{i}- Y, \A) \ar[d]^-{\pi'} \\
\bigoplus_{\varepsilon> 0} \gc^{\varepsilon}(X, \A) \ar[r]^-{v} & \gc^{0}(X, \A) \ar[r]^{\sim} & \gc(X- Y, \A) } \]
where the isomorphisms on the right come from \ref{pair}. We wish to show that $v$ is epi. Since $u$ is epi by the Hausdorff case, it suffices to show that
$\pi$ is epi, or, equivalently, that $\pi'$ is epi. This is indeed the case by \ref{MV}.\ \ $\Boxe$ \\
\newline
\hspace*{.3in} It is quite clear that using c-soft resolutions one can define compactly supported cohomology $H^{*}_{c}(X, \A)$ for any $\A \in \underline{Ab}(X)$. In particular, we get an extension $H^{0}_{c}(X, -)$ of $\gc(X, -)$ to all sheaves; this extension is still denoted by $\gc(X, -)$.


\begin{prop}
\label{losr}
Let $f: Y\rmap X$ be a continuous map. There is a functor $f_{\, !}: \underline{Ab}(Y) \rmap \underline{Ab}(X)$ with the following properties: \\
 \hspace*{.3in} \ps i\pd\ \ \ For any open $U \subset X$ and any $\B \in \underline{Ab}(Y)$, $\gc (U, f_{\, !}\B) = \gc (f^{-1}(U), \B)$.\\
 \hspace*{.3in} \ps ii\pd\ \ For any point $x\in X$ and any $\B \in \underline{Ab}(Y)$, $f_{\, !}(\B)_{x} = \gc(f^{-1}\ps x\pd, \B)$.\\
 \hspace*{.3in} \ps iii\pd\ $f_{\, !}$ is left exact and maps c-soft sheaves into c-soft sheaves.\\
 \hspace*{.3in} \ps iv\pd\ \ For any fibered product
\[ \xymatrix{
Z\times_{X}Y \ar[d]_-{q}\ar[r]^-{p} & Y \ar[d]^-{f}\\
Z\ar[r]^-{e} & X } \]
along an \'etale map $e$ and for any c-soft $\B\in \underline{Ab}(Y)$, there is a canonical isomorphism
\[ q_{\, !}p^{*}\B \cong e^{*}f_{\, !}\B \ .\]
(see \ref{pull} below for the case where $e$ is not \'etale).
\end{prop}

\emph{Proof:} Of course the proposition is well known in the Hausdorff case. For the more general case, recall first from \cite{Bredon} the correspondence for any {\it Hausdorff} space $Z$ between c-soft sheaves $\es$ on $Z$ and flabby cosheaves $\C$ on $Z$, given by:
\begin{equation}
\label{relA4}
 \gc(W, \es) = \C(W) 
\end{equation}
(natural with respect to the opens $W\subset Z).$ Given the cosheaf $\C$, the stalk of the corresponding sheaf $\es$ at a point $z\in Z$ is given by the exact sequence:
\begin{equation}
\label{relA5}
 0\rmap \C(Z- z)\rmap \C(Z)\rmap \es_{z}\rmap 0\ .
\end{equation}
\hspace*{.3in} We use this correspondence in the construction of $f_{\, !}$. (However, see remark \ref{rk} below for a description of $f_{\, !}$ which doesn't use this correspondence).\\
\hspace*{.3in} We discuss first the construction of $f_{\, !}$ on c-soft sheaves. Let $\B\in \underline{Ab}(Y)$ be c-soft. First, assume $X$ is Hausdorff. Let $\B$ be a c-soft sheaf on $Y$, and define a cosheaf $\C= c(\B)$ by $\C(U)= \gc(f^{-1}(U), \B)$. Note that $\C$ is indeed a flabby cosheaf, by \ref{MV}. By the correspondence (\ref{relA4}), there is a c-soft sheaf $\es$ on $X$, uniquely determined up to isomorphism by the identity $\gc(U, \es)= \C(U)$ for any open $U\subset X$. Thus, if $X$ is Hausdorff, we can define $f_{\, !}\B$ to be this sheaf $\es$.\\
\hspace*{.3in} In the general case, cover $X$ by Hausdorff opens $U_{i}$, and define in this way for each $i$ a c-soft sheaf $\es_{i}$ on $U_{i}$ by:
\begin{equation}
\label{relA6} \gc(V, \es_{i}) = \gc(f^{-1}(V), \B)\ .
\end{equation}
Then (again by the equivalence between sheaves and cosheaves) there is a canonical isomorphism $\theta_{i\, j}: \es_{j}|_{U_{i\, j}} \rmap \es_{i}|_{U_{i\, j}}$ satisfying the cocycle condition. Therefore  the sheaves $\es_{i}$ patch together to a sheaf $\es$ on $X$, uniquely determined up to isomorphism by the condition that $\es|_{U_{i}} = \es_{i}$ (by an isomorphism compatible with $\theta_{i\, j}$). Thus we can define $f_{\, !}\B$ to be $\es$.\\
\hspace*{.3in} We prove the properties $\ps i\pd -\ps iv\pd$ in the statement of the proposition for $\B\in \underline{Ab}(Y)$ c-soft. Property $\ps i\pd$ clearly holds for an open set $U$ contained in some $U_{i}$, by (\ref{relA6}). For general $U$, property $\ps i\pd$ then follows by the Mayer-Vietoris sequence. Next, identity (\ref{relA5}) yields for any point $x\in X$ an exact sequence:
\[ 0 \rmap \gc(Y-f^{-1}\ps x\pd, \ \B) \rmap \gc(Y, \B) \rmap f_{\, !}(\B)_{x} \rmap 0\ ,\]
and hence, by \ref{pair} the isomorphism $\ps ii\pd$. of the Proposition. Finally, $\ps iv\pd$ is clear from the local nature of the construction of $f_{\, !}$.\\ 
\hspace*{.3in} For general $\A\in \underline{Ab}(Y)$ we define $f_{\, !}(\A) \in \underline{Ab}(X)$ as the kernel of the map $f_{\, !}(\es^{0}) \rmap f_{\, !}(\es^{1})$ where $0\rmap \A\rmap \es^{0}\rmap \es^{1}\rmap . . . $ is a c-soft resolution of $\A$ (from the first part it follows that it is well defined up to isomorphisms). The properties $\ps i\pd  $ and $\ps ii\pd $  are now immediate consequences of the definition and of the previous case. Using \ref{qi} and $\ps ii\pd$ it easily follows that $f_{\, !}$ transforms acyclic complexes of c-soft sheaves on $\underline{Ab}(Y)$ into acyclic complexes on $\underline{Ab}(X)$. This immediately implies that $\f_{\, !}$ is left exact.\ \ $\Boxe$

\begin{num}
\label{rk}
{\bf Remark.}\emph{ We outline an alternative construction and proof of Proposition \ref{losr}, which does not use the correspondence between sheaves and cosheaves. This construction will be used in the proof of \ref{pull} below. We will assume that $\B$ is c-soft and $X$ is Hausdorff. (As in the proof of \ref{losr}, the construction of $f_{\, !}$ for general $X$ is then obtained by glueing the constructions over a cover by Hausdorff opens $U_{i}\subset X$.)\\
\hspace*{.3in} So, let $\B$ be a c-soft sheaf on $Y$. For any open set $V\subset Y$, denote by $\B_{V}$ the sheaf on $Y$ obtained by extending $\B|_{V}$ by zero. Thus $\B_{V}$ is evidently c-soft, and $\gc(Y, \B_{V}) = \gc(V, \B)$. Moreover, an inclusion $V\subset W$ induces an evident map $\B_{V}\hookrightarrow  \B_{W}$.\\
\hspace*{.3in}  Now let $Y= \bigcup W_{i}$ be a cover by Hausdorff open sets. This cover induces a long exact sequence:}
\[ .\ .\ .\ \rmap \bigoplus_{i_{0}< i_{1}} \B_{W_{i_{0}i_{1}}} \rmap \bigoplus_{i_{0}} \B_{W_{i_{0}}} \rmap \B \rmap 0\]
\emph{of c-soft sheaves on $Y$. By Corollary \ref{qi}, the functor $\gc(Y, -)$ applied to this long exact sequence again yields an exact sequence, and this is precisely the Mayer-Vietoris sequence of \ref{MV}. For each $i_{0}, . . . ,i_{n}$ let $f_{i_{0}, . . . ,i_{n}}: W_{i_{0}, . . . ,i_{n}}\rmap X$ be the restriction of $f$; this is a map between Hausdorff spaces, so we have $(f_{i_{0}, . . . ,i_{n}})_{\, !}(\B_{W_{i_{0}, . . . ,i_{n}}})$ defined as usual. Define $f_{\, !}(\B)$ as the cokernel fitting into a long exact sequence:}
\begin{equation}
\label{relA7}
 .\ .\ .\ \rmap \bigoplus_{i_{0}< i_{1}} (f_{i_{0}i_{1}})_{\, !}(\B_{W_{i_{0}i_{1}}}) \rmap \bigoplus_{i_{0}} (f_{i_{0}})_{\, !}(\B_{W_{i_{0}}}) \rmap f_{\, !}(\B) \rmap 0\ .
\end{equation}
\emph{\hspace*{.3in}For $ x\in X$, we have $(f_{i_{0}})_{\, !}(\B_{W_{i_{0}}})_{x}=$ $\gc (f^{-1}\ps x\pd \cap W_{i_{0}}; \B)$ by the Hausdorff case. So taking stalks of the long exact sequence in (\ref{relA7}) at $x$ and using the Mayer-Vietoris sequence \ref{MV} for the space $f^{-1}\ps x\pd$ we find $f_{\, !}(\B)_{x}=$ $\gc (f^{-1}\ps x\pd, \B)$ as in \ref{losr} 
$\ps ii\pd $. Property \ref{losr} $\ps i\pd$ is proved in a similar way (using \ref{pair}).\\
\hspace*{.3in} The functor $f_{\, !}$ can be extended to the derived category $D(Y)$ by taking a c-soft resolution $0\rmap \A\rmap \es^{0}\rmap \es^{1}\rmap . . . $ and defining $\er f_{\, !}(\A)$ as the complex $f_{\, !}(\es^{\cdot})$. Up to quasi-isomorphisms, this complex is well defined and does not depend on the resolution $\es^{\cdot}$, by \ref{pair}. (In this way, we obtaine in fact a well defined functor $\er f_{\, !}: D(Y) \rmap D(X)$ at the level of derived categories, which is sometimes simply denoted by $f_{\, !}$ again). In particular, $\cal{H}^{*}(\er \f_{\, !}(\A))$ gives in fact the right derived functors $R^{*}f_{\, !}$ of $f_{\, !}$.}
\end{num}

\begin{lem}
\label{pull}
For any pullback diagram:
\[ \xymatrix{
Z\times_{X}Y \ar[d]_-{q}\ar[r]^-{p} & Y \ar[d]^-{f}\\
Z\ar[r]^-{e} & X } \]
and any sheaf $\B$ on $Y$, there is a canonical quasi-isomorphism:
\[ (\er q_{\, !})p^{*}\B \simeq e^{*}(\er f_{\, !})\B \ .\]
\end{lem}

\emph{Proof:} Using Mayer-Vietoris for covers of $X$ and $Z$ by Hausdorff open sets, it suffices to consider the case where $X$ and $Z$ are both Hausdorff. Clearly it also suffices to prove the lemma in the special case where $\B$ is c-soft.\\
\hspace*{.3in} Let $Y= \bigcup W_i$ as in \ref{rk}, so that $f_{\, !}(\B)$ fits into a long exact sequence (\ref{relA7}) of c-soft sheaves on $X$. Applying the exact functor $e^*$ to this sequence and using the lemma in the Hausdorff case, one obtains a long exact sequence of the form:
\begin{equation} 
\label{cinci}
  .\ .\ .\ \rmap \bigoplus_{i_{0}< i_{1}} q_{\, !}p^*(\B|_{W_{i_0i_1}}) \rmap \bigoplus_{i_{0}} q_{\, !}p^*(\B|_{W_{i_0}}) \rmap e^*f_{\, !}(\B) \rmap 0\ .
\end{equation}
Now let $p^*(\B)\rmap \es^{\cdot}$ be a c-soft resolution over the pullback $Z\times_XY$. Then for any open $U\subset Y$, $\es_{p^{-1}(W)}^{\cdot}$ is a c-soft resolution of $p^*(\B_{W})$, so $q_{\, !}(\es_{p^{-1}(W)}^{\cdot})$ is a c-soft resolution of $q_{\, !}p^*(\B)$. The lemma now follows by comparing the sequence (\ref{cinci}) to the defining sequence
\[   .\ .\ .\ \rmap \bigoplus_{i_{0}< i_{1}} q_{\, !}(\es_{p^{-1}W_{i_0i_1}}) \rmap \bigoplus_{i_{0}} q_{\, !}(\es_{p^{-1}W_{i_0}}) \rmap q_{\, !}(\es) \rmap 0 \]
for $q_{\, !}(p^*(\B))\stackrel{def}{=} q_{\, !}(\es)$.  \ \ $\Boxe$

\begin{num}
\label{lset}
\emph{
{\bf $f_{\, !}$ on \'etale maps.} Let $f:Y \rmap X$ be an \'etale map, i.e. a local homeomorphism. It is well known that the pullback functor $f^{*}: \underline{Ab}(X) \rmap \underline{Ab}(Y)$ has an exact left-adjoint $f_{\, !}: \underline{Ab}(Y) \rmap \underline{Ab}(X)$, described on the stalks by $f_{\, !}(\B)_{x}= \oplus_{y\in f^{-1}\ps x\pd} \B_{y}$. This construction agrees with the one in \ref{losr}. In particular, for \'etale $f$, the counit of the adjunction defines a map: 
\[ \Sigma_{f}: f_{\, !}f^{*}(\A) \rmap \A\ ,\]
``summation along the fiber'', for any sheaf $\A$ on $X$.
}
\end{num}

\begin{num}
\label{lspr}
\emph{
{\bf $f_{\, !}$ on proper maps.} Define a map $f: Y\rmap X$ between (non-necessarily Hausdorff) spaces to be {\it proper} if:\\
\hspace*{.3in} (i) the diagonal $Y \rmap Y\times_{X} Y$ is closed.\\
\hspace*{.3in} (ii) for any Hausdorff open $U \subset X$ and any compact $K \subset U$, the set $f^{-1}(K)$ is compact.\\
It is easy to see that if $f$ is proper then $f_{\, !}= f_{*}$, as in the Hausdorff case. Furthermore, for any c-soft sheaf $\A$ on $X$, there is a natural map $\gc (X, \A) \rmap \gc(Y, f^{*}\A)$ defined by pullback, as in the Hausdorff case.
}
\end{num}

\begin{num}
\emph{
{\bf Remark:} Although this does not simplify matters, one could theoretically interpret some of the constructions and results of this Appendix as follows. First, observe that for {\it Hausdorff} groupoids, the results in Sections 1-6 of the paper can be based on the usual definition of $\Gamma_c$ and are independent of the Appendix. Now, any non-separated manifold (or sufficiently nice space, cf. \ref{assumptions}) $X$ can be viewed as a trivial groupoid (\ref{ex1}.1), which is Morita equivalent to the Hausdorff \'etale groupoid $\G$ defined from an open cover $\{ U_i\}$ of $X$ by Hausdorff open sets, by taking $\nG{0}$ to be the disjoint sum of the $U_i$, and $\nG{1}= \nG{0}\times_X\nG{0}$.
}
\end{num}

Marius Crainic,\\
\hspace*{.2in}Utrecht University, Department of Mathematics,\\
\hspace*{.2in}P.O.Box:80.010,3508 TA Utrecht, The Netherlands,\\ 
\hspace*{.2in}e-mail: crainic@math.ruu.nl\\
\newline

Ieke Moerdijk,\\
\hspace*{.2in}Utrecht University, Department of Mathematics,\\ 
\hspace*{.2in}P.O.Box:80.010,3508 TA Utrecht, The Netherlands,\\ 
\hspace*{.2in}e-mail: moerdijk@math.ruu.nl


\begin{thebibliography}{xxxx}

\bibitem{alv}{\sc Jes\'us A. Alvarez L\'opez, }{\it A decomposition theorem for the spectral sequence of Lie foliations, }{Trans. of the A.M.S., Volume 329, Number 1, January 1992, 173- 184}
\bibitem{SGA}{\sc M. Artin, A. Grothendieck {\rm and} J.L. Verdier, }{\it Th\'eorie de topos et cohomologie \'etale des sch\'emas, SGA4, tome 2, Springer LNM }{\bf 270}{, 1972}
\bibitem{BoTu}{\sc R. Bott {\rm and} L.W. Tu, }{\it Differential forms in algebraic topology, }{Graduate Texts in Mathematics 82, Springer }{(1982/1995)}
\bibitem{Bo}{\sc R. Bott, }{\it Lectures on characteristic classes and foliations, }{Springer LNM }{\bf 279}{, 1-94}
\bibitem{Borel}{\sc Borel, }{\it Intersection cohomology, }{Birkhauser, 1984}
\bibitem{Bredon}{\sc Bredon, }{\it Sheaf Theory, }{McGraw-Hill, 1967} 
\bibitem{BrNi}{\sc J.-L. Brylinski {\rm and} V. Nistor, }{\it Cyclic Cohomology of Etale Groupoids, }{K-theory }{\bf 8 }{(1994), }{341-365}
\bibitem{BufLor}{\sc J.P. Buffet {\rm and} J.C. Lor, }{\it Une construction d'un universel pour une classe assez large de $\Gamma$-structures, }{C.R. Acad. Sci. Paris Ser. A-B }{\bf 270 }{(1970), A640- A642}
\bibitem{Bu}{\sc D. Burghelea, }{\it The cyclic homology of the group rings, }{Comm. Math. Helv. }{\bf 60 }{(1985), }{354-365}
\bibitem{CoOp}{\sc A. Connes, }{\it A survey of foliations and operator algebras, }{Proc. Sympos. Pure Math., AMS Providence, }{\bf 32}{ (1982)}{, 521-628}
\bibitem{Co1}{\sc A. Connes, }{\it Cohomologie cyclique et foncteur $Ext^{n}$, }{CRAS Paris }{\bf 296 }{(1983), }{953-958}
\bibitem{Co2}{\sc A. Connes, }{\it Noncommutative differential geometry Ch. II: de Rham homology and non commutative algebra, }{Publ. Math. IHES }{\bf 62 }{(1985),}{94-144}
\bibitem{Co3}{\sc A. Connes, }{\it Noncommutative Geometry, }{Academic Press }{(1994)}
\bibitem{Cra}{\sc M. Crainic, }{\it Cyclic homology of \'etale groupoids; The general case, }{K Theory, to appear}
\bibitem{FeTs}{\sc B.P. Feigin {\rm and} B.L. Tsygan, }{\it Additive K-theory, }{LNM }{\bf 1289 }{(1987) }{ }
\bibitem{God}{\sc R. Godement, }{\it Topologie Alg\'ebrique et Th\'eorie des Faisceaux, }{Hermann Paris }{1958}
\bibitem{Gr}{\sc A. Grothendieck, }{\it Produit tensoriel topologique et espaces nucleaires, }{Memoires AMS }{\bf 16 }{1955}
\bibitem{Ha1}{\sc A. Haefliger, }{\it Homotopy and integrability, }{in: Manifolds, Amsterdam, 1970, }{SLN }{\bf 192}{ (1971), 133- 163}
\bibitem{difcoh}{\sc A. Haefliger, }{\it Differentiable cohomology, }{"Differential Topologi (CIME, Varenna 1976) Ligouri, Naples", }{1979, }{19-70}{}
\bibitem{Ha2}{\sc A. Haefliger, }{\it Groupo\"\i des d'holonomie et espaces classifiants, }{Ast\'erisque }{\bf 116 }{(1984), }{70-97}
\bibitem{man}{\sc A. Haefliger, }{\it Cohomology theory for \'etale groupoids, }{unpublished manuscript }{(1992)}{}
\bibitem{HiSk}{\sc M. Hilsum {\rm and} G. Skandalis, }{\it Morphismes K-orientes d'espaces de feuilles et functorialite en theorie de Kasparov, }{Ann. Scient. ENS }{\bf 20}{ (1987), 325- 390}
\bibitem{Iv}{\sc B. Iversen, }{\it Sheaf Theory, }{Springer-Verlag, 1986}{}
\bibitem{Kar1}{\sc M. Karoubi, }{\it Homologie cyclique des groupes et des algebres, }{CRAS Paris }{\bf 297 }{(1983), }{381-384}
\bibitem{Kar2}{\sc M. Karoubi, }{\it Homologie cyclique et $K$-th\'eorie, }{Ast\'erisque }{\bf 149 }{(1987)}
\bibitem{Kassel}{\sc C. Kassel, }{\it Cyclic homology, comodules and mixed complexes, }{J. of Algebra }{\bf 107}{ (1987), 195- 216}
\bibitem{Kaw}{\sc T. Kawasaki, }{\it The index of elliptic operators over V-manifolds, }{Nagoya Math. Journal }{\bf 84}{ (1981), 135-157}
\bibitem{MaL}{\sc S. MacLane, }{\it Homology, Berlin, Heidelberg, New York: Springer Verlag, 1963}
\bibitem{MP}{\sc I. Moerdijk {\rm and} D. A. Pronk, }{\it Orbifolds, Sheaves and Groupoids, }{K-Theory }{\bf 12 }{(1997), }{3-21}
\bibitem{Fourier}{\sc I. Moerdijk, }{\it Classifying topos and foliations, }{Ann. Inst. Fourier. Grenoble }{\bf 41 }{(1991), }{189-209}
\bibitem{emb}{\sc I. Moerdijk, }{\it Proof of a conjecture of A. Haefliger, }{Topology }{\bf 37 }{(1998), 735-741}
\bibitem{Mr}{\sc J. Mrcun, }{\it Stability and invariants of Hilsum-Skandalis maps, }{Ph.D.-thesis, Utrecht Univ. }{(1996)}
\bibitem{Mrr}{\sc J. Mrcun, }{\it Functoriality of the bimodule associated to a Hilsum-Skandalis map}{, preprint, 1996}
\bibitem{Ni1}{\sc V. Nistor, }{\it Group cohomology and the cyclic cohomology of crossed products, }{Invent. Math. }{\bf 99 }{(1990), }{411-424}
\bibitem{Segal}{\sc G.Segal, }{\it Classifying spaces and spectral sequences, }{Publ. Math. IHES }{\bf 34 }{(1968), }{105-112}
\bibitem{Sp}{\sc E.H. Spanier, }{\it Algebraic Topology, }{Mc. Graw Hill, New York, etc, 1966}
\bibitem{We1}{\sc C. Weibel, }{\it An introduction to homological algebra, }{Cambridge Studies in Advanced Mathematics 38 }{(1994)}
\bibitem{Wi}{\sc H. Winkelnkemper, }{\it The graph of a foliation, }{Ann. Global Anal. Geom. }{\bf 1 }{(1983), }{51-75}

\end{thebibliography}
\end{document}